\documentclass[a4paper]{article}

\usepackage{etex}

\usepackage{amssymb}
\usepackage{amsmath}
\usepackage{mathtools}
\usepackage[amsmath,thmmarks,hyperref]{ntheorem}
\usepackage{enumerate}
\usepackage{dsfont}
\usepackage{color}

\usepackage{tikz}
\usepackage{pgfplots}
	\pgfplotsset{compat=1.9}

\usepackage{booktabs}

\usepackage{geometry}
\usepackage[final,notref]{showkeys}	
\usepackage[bookmarksnumbered=true,hidelinks]{hyperref}
\usepackage{caption}

\theoremnumbering{arabic}
\theoremstyle{break}
\RequirePackage{latexsym}
\theorembodyfont{\itshape}
\theoremsymbol{}
\theoremheaderfont{\normalfont\bfseries}
\theoremseparator{}
	\newtheorem{theorem}{Theorem}[section]
	\newtheorem{lemma}[theorem]{Lemma}

\theorembodyfont{\upshape}
\theoremsymbol{\ensuremath{\clubsuit}}
	\newtheorem{remark}[theorem]{Remark}
	
\theorembodyfont{\upshape}
\theoremsymbol{}
	\newtheorem{example}[theorem]{Example}
	
\theoremstyle{nonumberplain}
\theoremheaderfont{\normalfont\bfseries}
\theoremseparator{:}
\theorembodyfont{\normalfont}
\theoremsymbol{\ensuremath{\square}}
\RequirePackage{amssymb}

        \newtheorem{proof}{Proof}
		\qedsymbol{\ensuremath{\square}}

\numberwithin{equation}{section}

\newcommand{\KAO}{\mathcal{O}}
\newcommand{\KAK}{\mathcal{K}}
\newcommand{\eps}{\varepsilon}
\newcommand{\tnorm}[1]{\left|\!\!\;\left|\!\!\;\left| {#1} \right|\!\!\;\right|\!\!\;\right|}
\newcommand{\norm}[1]{\left\| {#1} \right\|}
\newcommand{\bignorm}[1]{\big\| {#1} \big\|}
\newcommand{\snorm}[1]{\left| {#1} \right|}

\title{Analysis of Galerkin and streamline-diffusion FEMs on piecewise equidistant
		meshes for turning point problems exhibiting an interior layer}
\author{Simon Becher\footnote{Institute of Numerical Mathematics, 
Technical University of Dresden, Dresden 01062, Germany.
\mbox{e-mail:} Simon.Becher@tu-dresden.de}}
\date{}

\begin{document}

	\maketitle
	
	\begin{abstract}
		We consider singularly perturbed boundary value problems with a simple interior
		turning point whose solutions exhibit an interior layer. These problems are discretised
		using higher order finite elements on layer-adapted piecewise equidistant meshes
		proposed by Sun and Stynes. We also study the streamline-diffusion finite element
		method (SDFEM) for such problems. For these methods error estimates uniform
		with respect to $\varepsilon$ are proven in the energy norm and in the stronger
		SDFEM-norm, respectively. Numerical experiments confirm the theoretical findings.
	\end{abstract}
	
	

	\section{Introduction}
	
	We consider singularly perturbed boundary value problems of the form
	\begin{subequations}
	\label{prob:intLayer}
	\begin{equation}
		\begin{aligned}
			-\varepsilon u''(x) + a(x) u'(x) + c(x) u(x) &= f(x) \qquad \text{in } (-1,1), \\
			u(-1) = \nu_{-1}, \quad u(1) &= \nu_1,
		\end{aligned}
	\end{equation}
	with a small parameter $0 < \varepsilon \ll 1$ and sufficiently smooth data $a, c, f$
	satisfying
	\begin{equation}
		a(x) = - (x-x_0) b(x), \qquad b(x) > 0, \qquad c(x) \geq 0, \qquad c(x_0) > 0
	\end{equation}
	\end{subequations}
	for a point $x_0 \in (-1,1)$. The simple zero $x_0$ of $a$ is an attractive
	simple turning point of the problem. Thus, the solution of~\eqref{prob:intLayer} exhibits an
	interior layer of ``cusp''-type~\cite{SS94} at $x_0$.
	
	In the literature (see e.g.~\cite{BHK84}, \cite[p.~95]{Lis01}, \cite[Lemma~2.3]{SS94})
	bounds for such interior layers are well known. We have
	\begin{gather}
		\label{ieq:innerLayerBounds}
		\left| u^{(i)}(x) \right| \leq C \left( 1 + \left(\eps^{1/2} + |x-x_0| \right)^{\lambda-i}\right)
	\end{gather}
	where the parameter $\lambda$ satisfies $0 < \lambda < \bar{\lambda} := c(x_0)/|a'(x_0)|$.
	Note that the estimate also holds for $\lambda = \bar{\lambda}$, if $\bar{\lambda}$ is not
	an integer. Otherwise there is an additional logarithmic factor, see references above.
	In the following we assume $x_0 = 0$ for convenience.
	
	The quest for uniform error estimates for singularly perturbed problems has concerned researchers
	for many years. One of the common strategies is the use of layer-adapted meshes to treat the
	occurring boundary and interior layers. In particular, meshes for layers of exponential type have
	been examined, see e.g.~\cite{Lin10} where various problems, numerical methods, and meshes
	are presented. Popular examples are, due to their simplicity, the piecewise equidistant Shishkin
	meshes~\cite{Shi90,MOS96} which are fine only in the layer region. Unfortunately, the layers of
	``cusp''-type~\eqref{ieq:innerLayerBounds} do not fade away that quickly and, thus, local
	refinements do not suffice to capture the layer. Therefore, Sun and Stynes~\cite[Section~5.1]{SS94}
	generalise the standard Shishkin approach and propose a mesh consisting of $\KAO(\ln N)$
	equidistant parts to analyse linear finite elements. Moreover, in~\cite{Lis01}
	Liseikin uses graded meshes adapted to~\eqref{ieq:innerLayerBounds} to prove the $\eps$-uniform
	first order convergence of an upwind scheme in the discrete maximum norm.
	
	For problems of the form~\eqref{prob:intLayer} the finite element method is analysed
	in~\cite{Bec16ArXivLis} on the graded meshes of Liseikin. Using related techniques we
	shall extend the results of Sun and Stynes~\cite{SS94} by studying finite elements of
	order $k \geq 1$ on piecewise uniform meshes with slightly modified parameters, see
	Section~\ref{sec:energySunStynes}. In particular, we prove $\eps$-uniform error estimates
	of the form $\left(N^{-1} \ln N\right)^k$ in a weighted energy norm.
	
	In numerical experiments non-physical oscillations in the error can be observed. In order to
	damp such behaviour various stabilisation techniques have been
	proposed in recent years. We shall study the streamline-diffusion finite element method (SDFEM)
	first introduced by Hughes and Brooks~\cite{HB79}. In Section~\ref{sec:SDFEM} we
	prove an error estimate in the SDFEM-norm. Moreover, for linear elements a supercloseness result
	is given which allows to improve the bound for the $L^2$-norm error. As an example for the analysis
	in the context of Shishkin meshes we may refer to Stynes and Tobiska~\cite{ST08} who studied a
	two-dimensional convection-diffusion problem with exponential boundary layers for $Q_p$-elements.
	
	Some numerical results are given to illustrate the theoretical findings in Section~\ref{sec:numExp}.
	
	Notation: In this paper let $C$ denote a generic constant independent of $\eps$
	and the number of mesh points. Furthermore, for an interval $I$ we use the usual
	Sobolev spaces $H^1(I)$, $H_0^1(I)$, $W^{k,\infty}(I)$, and $L^2(I)$. The space
	of continuous functions on $I$ is written as $C(I)$. We denote by
	$\left(\cdot,\cdot\right)_I$ the usual $L^2(I)$ inner product and by 
	$\norm{\cdot}_I$ the $L^2(I)$-norm. Moreover, the supremum norm on $I$ is written as
	$\norm{\cdot}_{\infty,I}$ and the seminorm in $H^1(I)$ as $\snorm{\cdot}_{1,I}$.
	If $I = (-1,1)$, the index $I$ in inner products, norms, and seminorms will be
	omitted. Further notation will be introduced later at the beginning of
	the sections where it is needed.
	
	\section{FEM-analysis on arbitrary meshes}

	The following section is based on the paper of Sun and Stynes~\cite{SS94}. While their approach
	merely allows the analysis of linear finite elements, the subsequent results enable the
	analysis of finite elements of higher order. We will only consider homogeneous
	Dirichlet boundary conditions $\nu_{-1} = \nu_1 = 0$. This is no restriction
	since it can be easily ensured by modifying the right-hand side $f$. Without loss
	of generality (cf.~\cite[Lemma 2.1]{SS94}) we may assume that
	\begin{gather}
		\label{gamma}
		\left( c - \tfrac{1}{2}a'\right)(x) \geq \gamma > 0, \qquad \text{ for all } x \in [-1,1], \quad
			\eps \text{ sufficiently small.}
	\end{gather}
	
	For $v, w \in H_0^1((-1,1))$ we set
	\begin{gather*}
		B_\eps\!\left(v,w\right) := \left(\eps v',w'\right) + \left(a v', w\right) + \left(c v,w\right)\!.
	\end{gather*}
	Thanks to~\eqref{gamma} the bilinear form $B_\eps\!\left(\cdot,\cdot\right)$ is uniformly
	coercive over $H_0^1((-1,1)) \times H_0^1((-1,1))$ in terms of the weighted energy norm
	$\tnorm{\cdot}_\eps$ defined by
	\begin{gather*}
		\tnorm{v}_\eps := \left( \eps \snorm{v}_1^2 + \norm{v}^2\right)^{1/2}.
	\end{gather*} 
	
	The weak formulation of~\eqref{prob:intLayer} with $\nu_{-1} = \nu_1 = 0$ reads
	as follows: \medskip
	
	Find $u \in H_0^1((-1,1))$ such that
	\begin{gather}
		\label{weakprob:intLayer}
		B_\eps\!\left(u,v\right) = \left(f,v \right), \qquad \text{ for all } v \in H_0^1((-1,1)).
	\end{gather}
	
	Let $-1 = x_{-N} < \ldots < x_i < \ldots < x_N = 1$ define an arbitrary mesh on
	the interval $[-1,1]$. The mesh interval lengths are given by $h_i := x_i-x_{i-1}$.
	For $k \geq 1$ we denote by $P_k((x_a,x_b))$ the space of polynomial functions of
	maximal order $k$ over $(x_a,x_b)$. Furthermore, we define the trial and test
	space $V^N$ by
	\begin{gather*}
		V^N := \left\{ v \in C([-1,1]) : v|_{(x_{i-1},x_i)} \in P_k((x_{i-1},x_i))\, \forall i, \, v(-1) = v(1) = 0 \right\}\!.
	\end{gather*}
	The discrete problem is given by: \medskip
	
	Find $u_N \in V^N$ such that
	\begin{gather}
		\label{dprob:intLayer}
		B_\eps\!\left(u_N,v_N\right) = \left(f,v_N \right), \qquad \text{ for all } v_N \in V^N.
	\end{gather}
	
	Let $u_I$ denote the standard Lagrangian interpolation into $V^N$, using the
	mesh points and $k-1$ (arbitrary) inner interpolation points per interval.
	For example uniform or Gau{\ss}-Lobatto points could be chosen.
	
	Assuming $u \in W^{k+1,\infty}((x_{i-1},x_i))$,
	the standard interpolation theory leads to the error estimates: 
	For all $j = 0, \ldots, k+1$
	\begin{gather}
		\bignorm{(u-u_I)^{(j)}}_{\infty,(x_{i-1},x_i)}
			\leq C h_i^{k+1-j} \bignorm{u^{(k+1)}}_{\infty,(x_{i-1},x_i)} \label{standInt} \\
		\intertext{and}
		\bignorm{(u-u_I)^{(j)}}_{\infty,(x_{i-1},x_i)}
			\leq C \bigl\| u^{(j)} \bigr\|_{\infty,(x_{i-1},x_i)} \label{standIntInfty}
	\end{gather}
	where $C$ depends on the choice of the inner interpolation points.
	Furthermore, for all $v_N \in V^N$ the inverse inequality
	\begin{gather}
		\label{invIneq}
		\norm{v_N'}_{\infty,(x_{i-1},x_i)} \leq C h_i^{-1} \norm{v_N}_{\infty,(x_{i-1},x_i)}
	\end{gather}
	holds.
	
	In order to estimate the error of the finite element solution we use
	the splitting
	\begin{gather}
		\label{eq:splitting}
		u-u_N = (u-u_I) + (u_I-u_N).
	\end{gather}
	The next lemma shows that the energy norm of the second term can be estimated by
	knowing some interpolation error bounds only. The given approach works for finite elements
	of arbitrary order $k\geq 1$.
	\begin{lemma}
		\label{le:energyuIuN_Pk}
		Let $u$ be the solution of~\eqref{prob:intLayer} and $u_N$ the solution of~\eqref{dprob:intLayer}
		on an arbitrary mesh. Then we have
		\begin{gather*}
			\tnorm{u_I - u_N}_\eps \leq C \Big(\tnorm{u_I-u}_\eps + \norm{x (u_I-u)'}\Big).
		\end{gather*}
	\end{lemma}
	\begin{proof}
		Using the coercivity of $B_\eps\!\left(\cdot,\cdot\right)$ and orthogonality, we obtain
		\begin{gather}
			\label{ieq:coerBeps}
			\min\{\gamma,1\} \tnorm{u_I-u_N}^2_\eps \leq B_\eps(u_I-u_N,u_I-u_N) = B_\eps(u_I-u,u_I-u_N).
		\end{gather}
		The Cauchy Schwarz inequality yields
		\begin{align*}
			& B_\eps(u_I-u,u_I-u_N) \\
			& \qquad = \eps \bigl( (u_I-u)',(u_I-u_N)'\bigr)
				+ \bigl(a(u_I-u)',u_I-u_N\bigr) + \bigl(c(u_I-u),u_I-u_N\bigr) \\
			& \qquad \leq \sqrt{\eps} \snorm{u_I-u}_1 \sqrt{\eps}\snorm{u_I-u_N}_1
				+ \norm{a(u_I-u)'} \norm{u_I-u_N}
				+ \norm{c}_\infty \norm{u_I-u} \norm{u_I-u_N}.
		\end{align*}
		Additionally, using the fact that $a(x) = -xb(x)$, we get
		\begin{align*}
			B_\eps(u_I-u,u_I-u_N)
				& \leq \Bigl(\sqrt{\eps}\snorm{u_I-u}_1 + \norm{b}_\infty \norm{x(u_I-u)'}
					+ \norm{c}_\infty \norm{u_I-u}\Bigr) \tnorm{u_I-u_N}_\eps \\
				& \leq C \Bigl( \tnorm{u_I-u}_\eps + \norm{x(u_I-u)'}\Bigr) \tnorm{u_I-u_N}_\eps.
		\end{align*}
		Combining this and~\eqref{ieq:coerBeps} completes the proof.
	\end{proof}
	
	\begin{remark}
		For linear finite elements Sun and Stynes \cite[Lemma~5.2]{SS94} proved an estimate of the form
		\begin{gather*}
			\tnorm{u_I-u_N}_\eps \leq C \left(\norm{u-u_I}^{1/2} + \max_i{h_i^2}\right)\!,
		\end{gather*}
		see also~\cite[Lemma~3.7]{Bec16ArXivLis}. Aside from the fact that their argument works for
		linear elements only, such an estimate would not enable optimal estimates for finite elements
		of higher order.
	\end{remark}
	
	\begin{remark}
		In the setting of Lemma~\ref{le:energyuIuN_Pk} we also have
		\begin{gather*}
			\tnorm{u_I - u_N}_\eps \leq C \tnorm{u_I-u}_\epsilon
				+ C \left(\sum_{i=-N+1}^N h_i^{-2} \norm{x (u_I-u)}_{(x_{i-1},x_i)}^2\right)^{\!1/2}
		\end{gather*}
		which is proven in~\cite[Lemma~3.1]{Bec16ArXivLis}.
	\end{remark}
	
	\section{The piecewise equidistant meshes of Sun and Stynes}
	\label{sec:energySunStynes}
	
	This section is devoted to the study of the piecewise equidistant meshes proposed by
	Sun and Stynes in~\cite[Section~5.1]{SS94}. They generalise the standard
	approach of Shishkin and introduce a mesh that consists of $\KAO(\ln N)$ equidistant
	parts. Because of symmetry, the mesh will be described for $x \geq 0$ only.
	In order to enable the analysis of finite elements of order $k \geq 1$, we slightly
	modify the mesh parameters.
	
	In the following we assume that $\lambda$ from~\eqref{ieq:innerLayerBounds} lies in $(0,k+1)$
	which is the most difficult case. Otherwise all crucial derivatives of the solution could be
	bounded by a generic constant independent of $\eps$ and consequently optimal order $\eps$-uniform
	estimates could be proven with standard methods on uniform meshes.
	
	For $\eps \in (0,1]$ and given positive integer $N$ we set
	\begin{gather}
		\label{sigma}
		\sigma = \max\big\{\eps^{(1-\lambda/(k+1))/2}, N^{-(2k+1)} \big\}
	\end{gather}
	and
	\begin{gather}
		\label{KAK}
		\KAK = \left\lfloor 1 - \frac{\ln(\sigma)}{\ln(10)} \right\rfloor,
	\end{gather}
	where $\lfloor z \rfloor$ denotes the largest integer less or equal to $z$.
	
	The piecewise equidistant mesh is constructed as follows: The interval $(0,1]$ is partitioned
	into the $\KAK+1$ subintervals $(0,10^{-\KAK}], \, (10^{-\KAK},10^{-\KAK+1}], \ldots , (10^{-1},1]$.
	Then in a second step each of these subintervals is divided uniformly into $\lfloor N/(\KAK+1) \rfloor$
	parts. For simplicity we assume that $\lfloor N/(\KAK+1) \rfloor = N/(\KAK+1)$.
	Hence, by construction we have
	\begin{gather}
		\label{eq:hiKAK}
		h_i = (\KAK+1)10^{-\KAK}N^{-1}, \quad \text{for} \quad x_i \in (0,10^{-\KAK}]
	\end{gather}
	and
	\begin{gather}
		\label{eq:hiEll}
		h_i = 9(\KAK+1)10^{-l}N^{-1}, \quad \text{for} \quad x_i \in (10^{-l},10^{-l+1}]
			\quad \text{and} \quad l = 1,\ldots,\KAK.
	\end{gather}
	
	From~\eqref{sigma}, \eqref{KAK}, and the properties of the logarithm we see that
	\begin{gather*}
		\KAK + 1 \leq 2 - \frac{\ln(\sigma)}{\ln(10)}
			\leq 2 + \min\left\{\frac{1-\lambda/(k+1)}{2}\frac{|\ln(\eps)|}{\ln(10)},(2k+1)\frac{\ln(N)}{\ln(10)}\right\}.
	\end{gather*}
	For $N$ sufficiently large (dependent on $k$) this estimate yields $\KAK + 1 \leq N$. Furthermore,
	we obtain
	\begin{gather}
		\label{ieq:KAKtoLN}
		\KAK + 1 \leq C \ln N.
	\end{gather}
	From~\eqref{KAK} we have $\ln(10)-\ln(\sigma) \geq \KAK \ln(10) > -\ln(\sigma)$ which implies
	\begin{gather}
		\label{ieq:sigma}
		10^{-1} \sigma \leq 10^{-\KAK} < \sigma.
	\end{gather}
	
	The next lemma is a generalisation of~\cite[Lemma~5.3]{SS94} and provides some important
	basic results for the mesh intervals of the piecewise equidistant mesh.
	\begin{lemma}
		\label{le:sunStynesh^kBounds}
		Let $j=0,1$. The following inequalities hold
		\begin{align}
			h_i^{k+1-j}\left(x_{i-1}+\eps^{1/2}\right)^{\lambda-(k+1-j)}
				& \leq C \left((\KAK + 1)N^{-1}\right)^{k+1-j}, \quad \text{for} \quad x_i \in (10^{-\KAK},1], \\
			h_i^{k+1-j}\left(x_{i-1}+\eps^{1/2}\right)^{\lambda-(k+1-j)}
				& \leq C\left(i-1\right)^{-(k+1-j)}, \quad \text{for} \quad x_i \in (x_1,10^{-\KAK}]. \label{ieq:sunStynes_im1}
		\intertext{If $\sigma = \eps^{(1-\lambda/(k+1))/2}$, then}
			h_i^{k+1-j}\left(x_{i-1}+\eps^{1/2}\right)^{\lambda-(k+1-j)}
				& \leq C \left((\KAK + 1)N^{-1}\right)^{k+1-j}, \quad \text{for} \quad x_i \in (0,10^{-\KAK}].
		\end{align}
		In general, the mesh interval length can be bounded by
		\begin{gather*}
			h_i \leq (\KAK + 1)N^{-1}.
		\end{gather*}
		Furthermore, in the case of $\sigma = N^{-(2k+1)}$, we have
		\begin{gather}
			\label{ieq:sunStynesX_1}
			x_1 = h_1 \leq (\KAK+1) N^{-2(k+1)}.
		\end{gather}
	\end{lemma}
	\begin{proof}
		In order to prove the last estimate, let $\sigma = N^{-(2k+1)}$. Combining~\eqref{eq:hiKAK}
		with~\eqref{ieq:sigma} yields
		\begin{gather*}
			x_1 = h_1 = (\KAK + 1)10^{-\KAK}N^{-1} \leq (\KAK + 1) \sigma N^{-1} = (\KAK + 1) N^{-2(k+1)}.
		\end{gather*}
		The general estimates for $h_i$ follow from~\eqref{eq:hiKAK} and~\eqref{eq:hiEll}.
		
		We assume $\lambda - (k+1-j) < 0$ in the following. Otherwise
		$\left(x_{i-1}+\eps^{1/2}\right)^{\lambda-(k+1-j)} \leq C$ would
		allow to deduce the wanted bounds very easily.
		
		So, let $x_i \in (10^{-l},10^{-l+1}]$ for some $l \in \{1, \ldots, \KAK\}$.
		With~\eqref{eq:hiEll} we obtain
		\begin{align*}
			h_i^{k+1-j}\left(x_{i-1}+\eps^{1/2}\right)^{\lambda-(k+1-j)}
				& \leq \left(9 (\KAK + 1) 10^{-l} N^{-1} \right)^{k+1-j}
					\left(10^{-l}+\eps^{1/2}\right)^{\lambda-(k+1-j)} \\
				& \leq C \left((\KAK + 1) N^{-1}\right)^{k+1-j} \left(10^{-l}\right)^{k+1-j}
					\left(10^{-l}\right)^{\lambda-(k+1-j)} \\
				& = C \left((\KAK + 1) N^{-1} \right)^{k+1-j} 10^{-\lambda l}
					\leq C \left((\KAK + 1) N^{-1}\right)^{k+1-j}.
		\end{align*}
		For $x_i \in (x_1,10^{-\KAK}]$ the fact that the mesh is equidistant in this interval implies
		\begin{align*}
			h_i^{k+1-j}\left(x_{i-1}+\eps^{1/2}\right)^{\lambda-(k+1-j)}
				& = \left(x_{i-1} (i-1)^{-1}\right)^{k+1-j}
					\left(x_{i-1}+\eps^{1/2}\right)^{\lambda-(k+1-j)} \\
				& \leq (i-1)^{-(k+1-j)} x_{i-1}^\lambda \\
				& \leq (i-1)^{-(k+1-j)}.
		\end{align*}
		Finally, let $\sigma = \eps^{(1-\lambda/(k+1))/2}$ and $x_i \in (0,10^{-\KAK}]$. Using~\eqref{eq:hiKAK}
		and~\eqref{ieq:sigma} we get
		\begin{align*}
			h_i^{k+1-j}\left(x_{i-1}+\eps^{1/2}\right)^{\lambda-(k+1-j)}
				& \leq \left((\KAK + 1) \sigma N^{-1} \right)^{k+1-j} \eps^{(\lambda-(k+1-j))/2} \\
				& \leq C \left((\KAK + 1)N^{-1}\right)^{k+1-j} \eps^{\lambda j/(2(k+1))} \\
				& \leq C \left((\KAK + 1) N^{-1}\right)^{k+1-j}.
		\end{align*}
	\end{proof}
	
	The interpolation error on the layer-adapted piecewise equidistant mesh proposed
	by Sun and Stynes shall be bounded in the following lemma which is a generalisation
	of~\cite[Lemma~5.4]{SS94}. We also refer to~\cite[Lemma~3.3, Lemma~3.4]{Bec16ArXivLis}
	where a similar argumentation is used to estimate the interpolation error on
	special graded meshes.
	\begin{lemma}
		\label{le:IntErrorSunStynes}
		Let $u$ be the solution of problem~\eqref{prob:intLayer} and $u_I \in V^N$ be its interpolant
		on the piecewise equidistant mesh given by~\eqref{sigma} -- \eqref{eq:hiEll}. Then
		\begin{gather}
			\label{ieq:L2IntSunStynes}
			\norm{u-u_I} \leq C \left((\KAK + 1) N^{-1}\right)^{k+1}
		\end{gather}
		and
		\begin{gather}
			\label{ieq:energyIntSunStynes}
			\tnorm{u-u_I}_\eps + \norm{x(u-u_I)'} \leq C \left((\KAK + 1) N^{-1}\right)^{k}.
		\end{gather}
	\end{lemma}
	\begin{proof}
		Thanks to the symmetry of the problem, we shall consider only $x \in [0,1]$. Furthermore,
		we use $j \in \{0,1\}$ to switch between the $L^2$-norm term and the $\eps$-weighted
		$H^1$-seminorm term. The estimate for $\norm{x(u-u_I)'}$ is also covered by $j=1$.
		
		Let $x \in (x_{i-1},x_i)$ for some $i$, where $x_i \in (10^{-\KAK},1]$. Then
		\begin{align*}
			\big(\eps^{j/2}+jx\big) \left| (u-u_I)^{(j)}(x) \right|
				& \leq C \big(\eps^{j/2}+jx\big) h_i^{k+1-j} \bignorm{u^{(k+1)}}_{\infty,(x_{i-1},x_i)} \\
				& \leq C \big(\eps^{j/2}+j(x_{i-1}+h_i)\big) h_i^{k+1-j}
					\left(1 + \left( x_{i-1} + \eps^{1/2} \right)^{\lambda - (k+1)}\right) \\
				& \leq C \left((\KAK + 1) N^{-1}\right)^{k+1-j},
		\end{align*}
		where we used~\eqref{standInt}, \eqref{ieq:innerLayerBounds}, and Lemma~\ref{le:sunStynesh^kBounds}. Hence,
		\begin{gather*}
			\eps^j \int_{10^{-\KAK}}^1 \left(( u-u_I)^{(j)}(x)\right)^2 dx
				+ j^2 \int_{10^{-\KAK}}^1 \left(x( u-u_I)^{(j)}(x)\right)^2 dx
				\leq C \left((\KAK + 1) N^{-1}\right)^{2(k+1-j)}.
		\end{gather*}
		
		Now, let $x \in (x_{i-1},x_i)$, where $x_i\in(0,10^{-\KAK}]$.
		We consider two different cases for $\sigma$.
		
		First, if $\sigma = \eps^{(1-\lambda/(k+1))/2}$ then as above
		Lemma~\ref{le:sunStynesh^kBounds} yields
		\begin{align*}
			\big(\eps^{j/2}+jx\big)\left| (u-u_I)^{(j)}(x) \right|
				& \leq C \big(\eps^{j/2}+j(x_{i-1}+h_i)\big) h_i^{k+1-j}
					\left(1 + \left( x_{i-1} + \eps^{1/2} \right)^{\lambda - (k+1)}\right) \\
				& \leq C \left((\KAK + 1) N^{-1}\right)^{k+1-j}
		\end{align*}
		and therefore
		\begin{gather*}
			\eps^j \int_0^{10^{-\KAK}} \left( (u-u_I)^{(j)}(x)\right)^2 dx
				+ j^2 \int_0^{10^{-\KAK}} \left(x( u-u_I)^{(j)}(x)\right)^2 dx
				\leq C \left((\KAK + 1) N^{-1}\right)^{2(k+1-j)}.
		\end{gather*}
		
		Finally, let $\sigma = N^{-(2k+1)}$. The integral over $(0,x_1)$ is estimated directly.
		By~\eqref{standIntInfty}, \eqref{ieq:innerLayerBounds}, and Lemma~\ref{le:sunStynesh^kBounds},
		especially~\eqref{ieq:sunStynesX_1}, we have
		\begin{align*}
			\eps^j \int_0^{x_1} \left( (u-u_I)^{(j)}(x)\right)^2 dx
				& \leq C \eps^j \int_0^{x_1} \bigl\|u^{(j)}\bigr\|_{\infty,(0,x_1)}^2 dx \\
				& \leq C \eps^j \left(1+\eps^{(\lambda - j)/2}\right)^{\!2} x_1 \\
				& \leq C (\KAK+1) N^{-2(k+1)}
		\end{align*}
		and
		\begin{align*}
			\int_0^{x_1} \big(x(u-u_I)'(x)\big)^2 dx
				& \leq 2 \norm{x u'}_{\infty,(0,x_1)}^2 \int_0^{x_1} 1\, dx
					+ 2 \norm{u_I'}_{\infty,(0,x_1)}^2 \int_0^{x_1} x^2 dx \\
				& \leq 2 \norm{x u'}_{\infty,(0,x_1)}^2 x_1 
					+ C h_1^{-2} \norm{u_I}_{\infty,(0,x_1)}^2 x_1^3 \\
				& \leq C (\KAK+1) N^{-2(k+1)},
		\end{align*}
		where additionally the inverse inequality~\eqref{invIneq} and the stability
		$\norm{u_I}_\infty \leq C \norm{u}_\infty$ is applied.
		Now, let $x \in (x_{i-1},x_i) \subseteq (x_1,10^{-\KAK}]$. Then
		\begin{align*}
			\big(\eps^{j/2}+jx\big) \left| (u-u_I)^{(j)}(x) \right|
				& \leq C \big(\eps^{j/2}+j(x_{i-1}+h_i)\big) h_i^{k+1-j}
					\left(1 + \left( x_{i-1} + \eps^{1/2} \right)^{\lambda - (k+1)}\right) \\
				& \leq C \left(\big(\eps^{j/2}+j 10^{-\KAK}\big) h_i + \left(i-1\right)^{-(k+1-j)}\right),
		\end{align*}
		for $i=2, \ldots,N/(\KAK+1)$, where Lemma~\ref{le:sunStynesh^kBounds},
		especially~\eqref{ieq:sunStynes_im1}, is used. Hence,
		with~\eqref{ieq:sunStynesX_1}
		\begin{align*}
			& \eps^j \int_{x_1}^{10^{-\KAK}} \left( (u-u_I)^{(j)}(x)\right)^2 dx
				+ j^2 \int_{x_1}^{10^{-\KAK}} \left(x( u-u_I)^{(j)}(x)\right)^2 dx \\
			& \qquad \qquad \leq C \sum_{i=2}^{N/(\KAK+1)} h_i
				\left(\big(\eps^{j/2}+j 10^{-\KAK}\big) h_i + \left(i-1\right)^{-(k+1-j)}\right)^2 \\
			& \qquad \qquad \leq C x_1 \sum_{i=2}^{N/(\KAK+1)} \left(h_i + \left(i-1\right)^{-2}\right) \\
			& \qquad \qquad \leq C x_1 \left(10^{-\KAK} + \tfrac{\pi^2}{6}\right) \\
			& \qquad \qquad \leq C (\KAK+1) N^{-2(k+1)}.
		\end{align*}
		
		Combining the above estimates for $j=0$ and using symmetry on $[-1,0]$ we get~\eqref{ieq:L2IntSunStynes}.
		This estimate together with the above estimates for $j=1$ immediate gives~\eqref{ieq:energyIntSunStynes}.
	\end{proof}
	
	Now, we are able to prove the $\eps$-uniform error estimate in the energy norm
	for finite elements of order $k\geq 1$ on the piecewise equidistant mesh of Sun and Stynes.
	\begin{theorem}
		\label{th:energyErrorSunStynes}
		Let $u$ be the solution of~\eqref{prob:intLayer} and $u_N$ the solution
		of~\eqref{dprob:intLayer} on the piecewise equidistant mesh given
		by~\eqref{sigma} -- \eqref{eq:hiEll}. Then we have
		\begin{gather*}
			\tnorm{u-u_N}_\eps \leq C \left((\KAK + 1) N^{-1}\right)^k \leq C \left(N^{-1}\ln N\right)^k.
		\end{gather*}
	\end{theorem}
	\begin{proof}
		The bound in the energy norm follows easily using the splitting~\eqref{eq:splitting},
		the triangle inequality, Lemma~\ref{le:energyuIuN_Pk}, and~\eqref{ieq:energyIntSunStynes}.
		The second inequality is an immediately consequence of~\eqref{ieq:KAKtoLN}.
	\end{proof}
	
	\begin{remark}
		Under certain assumptions it was proven for special graded meshes proposed by Liseikin that
		\begin{gather*}
			\tnorm{u-u_N}_\eps \leq C N^{-k},
		\end{gather*}
		see~\cite[Theorem~3.5]{Bec16ArXivLis}. Thus, these meshes seem to be optimal in the sense that
		no additional logarithmic factor appears in the error estimate. However, the constant
		may depend on a mesh parameter $\alpha \in (0,\lambda]$.
		
		Furthermore, note that Liseikin's meshes are not well-defined for $\lambda = 0$ whereas
		the construction of the Sun and Stynes meshes works in this case as well. Therefore,
		the latter meshes can also be used to handle certain power-type layers caused by simple
		boundary turning points, for details see~\cite{Bec17ArXivGen}.
		
		The two types of meshes have been compared numerically in~\cite{Bec16MAIS}.
	\end{remark}
	
	\section{SDFEM-analysis on arbitrary and piecewise equidistant meshes}
	\label{sec:SDFEM}
	
	In this section the streamline-diffusion finite element method (SDFEM) is studied.
	For convenience we use the shorter notation $I_i:=(x_{i-1},x_i)$ for all $i = -N+1,\ldots,N$.
	In order to increase the stability some extra terms are added to the weak formulation. We set
	\begin{gather*}
		B_{\mathrm{SD}}\!\left(v,w\right) := B_{\eps}\!\left(v,w\right) + B_{\mathrm{Stab}}\!\left(v,w\right)
	\end{gather*}
	where
	\begin{gather*}
		B_{\mathrm{Stab}}\!\left(v,w\right)
			:= \sum_{i=-N+1}^N \delta_i \int_{x_{i-1}}^{x_i}
				\big(-\varepsilon v'' + a v' + c v\big)(x)\big(aw'\big)(x) dx
	\end{gather*}
	and
	\begin{gather*}
		f_{\mathrm{Stab}}\!\left(v\right)
			:= \sum_{i=-N+1}^N \delta_i \int_{x_{i-1}}^{x_i}
				f(x) \big(av'\big)(x) dx
	\end{gather*}
	with stabilisation parameters $\delta_i \geq 0$ to be defined later.
	Now, the discrete problem is given by: \medskip
	
	Find $u_N \in V^N$ such that
	\begin{gather}
		\label{dprob:SDFEM}
		B_{\mathrm{SD}}\!\left(u_N,v_N\right) = \left(f,v_N\right) + f_{\mathrm{Stab}}(v_N),
			\qquad \text{ for all } v_N \in V^N.
	\end{gather}
	Note that the method is consistent, i.e. for $u\in H^2((-1,1))$
	of~\eqref{weakprob:intLayer} we have
	\begin{gather*}
		B_{\mathrm{SD}}\!\left(u,v_N\right) = \left(f,v_N\right) + f_{\mathrm{Stab}}(v_N),
			\qquad \text{ for all } v_N \in V^N.
	\end{gather*}
	
	For our analysis we define the SDFEM-norm by
	\begin{gather*}
		\tnorm{v}_{\mathrm{SD}} := \left( \eps \snorm{v}_1^2 + \norm{v}^2
			+ \sum_{i=-N+1}^N \bignorm{\sqrt{\delta_i}av'}_{I_i}^2 \right)^{\!1/2}.
	\end{gather*}
	Because of the additional terms this norm is stronger than the energy norm.
	
	The following inverse inequality holds
	\begin{gather*}
		\norm{v_N''}_{I_i} \leq  c_{\mathrm{inv}}h_i^{-1} \norm{v_N'}_{I_i}\!, \qquad \text{ for all } v_N \in V^N
	\end{gather*}
	with a constant $c_\mathrm{inv}$ independent of $i$ and $h_i$.
	Thus, imposing the requirement
	\begin{subequations}
	\label{ass:delta}
	\begin{gather}
		0 \leq \delta_i
			\leq \frac{1}{2} \min \left\{\frac{h_i^2}{\eps c_{\mathrm{inv}}^2}, \frac{\gamma}{\norm{c}_{\infty}^2}\right\}
	\end{gather}
	or the assumption
	\begin{gather}
		0 \leq \delta_i
			\leq \frac{\gamma}{2 \norm{c}_{\infty}^2}
	\end{gather}
	\end{subequations}
	for linear finite elements, respectively, we obtain analogously to~\cite[p.~86]{RST08}
	\begin{gather}
		\label{ieq:coerSDFEM}
		B_{\mathrm{SD}}\!\left(v_N,v_N\right) \geq \frac{1}{2} \min\{\gamma,1\} \tnorm{v_N}_{\mathrm{SD}}^2\!,
			\qquad \text{ for all } v_N \in V^N.
	\end{gather}
	
	\subsection{Higher order finite elements}
	Using the splitting~\eqref{eq:splitting} our analysis starts with the following lemma.
	\begin{lemma}
		\label{le:SDFEMuIuN_Pk}
		Let $u$ be the solution of~\eqref{prob:intLayer}, $u_N$ the solution of~\eqref{dprob:SDFEM},
		and $u_I$ the interpolant of $u$ on an arbitrary mesh. Furthermore, choose $\delta_i$ such
		that~\eqref{ass:delta} is satisfied. Then we have
		\begin{align*}
			\tnorm{u_I - u_N}_{\mathrm{SD}}
				&\leq C \Big(1+\max_i \sqrt{\delta_i} \Big)
					\Big(\tnorm{u_I-u}_\eps + \norm{x (u_I-u)'} \Big) \\
				& \qquad + C \left(\sum_{i=-N+1}^N
					\min\left\{ \bignorm{\sqrt{\delta_i}\,\eps (u_I-u)''}_{I_i}^2,
						\norm{\delta_i \sqrt{\eps}\, x(u_I-u)''}_{I_i}^2\right\} \right)^{\!1/2}.
		\end{align*}
	\end{lemma}
	\begin{proof}
		By~\eqref{ieq:coerSDFEM}
		and due to orthogonality which is implied by consistency, we have
		\begin{gather*}
		\begin{aligned}
			\frac{1}{2} \min\{\gamma,1\} \tnorm{u_I-u_N}_{\mathrm{SD}}^2
				& \leq B_{\mathrm{SD}}\!\left(u_I-u_N,u_I-u_N\right)
					= B_{\mathrm{SD}}\!\left(u_I-u,u_I-u_N\right) \\ 
				& = B_\eps\!\left(u_I-u,u_I-u_N\right)
					+ B_{\mathrm{Stab}}\!\left(u_I-u,u_I-u_N\right).
		\end{aligned}
		\end{gather*}
		As in the proof of Lemma~\ref{le:energyuIuN_Pk} we obtain for the first term
		\begin{gather*}
			B_\eps(u_I-u,u_I-u_N)
				\leq C \Bigl( \tnorm{u_I-u}_\eps + \norm{x(u_I-u)'}\Bigr) \tnorm{u_I-u_N}_\eps.
		\end{gather*}
		
		It remains to estimate the second term
		\begin{gather}
		\label{ieq:Bstab}
		\begin{aligned}
			& B_{\mathrm{Stab}}\!\left(u_I-u,u_I-u_N\right) \\
			& \qquad = \sum_i \delta_i \int_{I_i}
				\big(-\varepsilon (u_I-u)'' + a (u_I-u)' + c (u_I-u)\big)(x)
				\big(a(u_I-u_N)'\big)(x) \,dx.
		\end{aligned}
		\end{gather}
		Applying Cauchy Schwarz' inequality and using the fact that $a(x)=-xb(x)$
		we gain for the last two summands in~\eqref{ieq:Bstab}
		\begin{gather}
		\label{ieq:BstabACterms}
		\begin{aligned}
			& \left|\sum_i \delta_i \int_{I_i}
				\big(a (u_I-u)' + c (u_I-u)\big)(x)
				\big(a(u_I-u_N)'\big)(x) \,dx \right|\\
			& \qquad \leq \sum_i \sqrt{\delta_i}
				\left(\norm{a(u_I-u)'}_{I_i}
					+ C\norm{u_I-u}_{I_i}\right)
					\bignorm{\sqrt{\delta_i} a(u_I-u_N)'}_{I_i} \\
			& \qquad \leq C \max_i \sqrt{\delta_i} \Big(\norm{x(u_I-u)'} + \norm{u_I-u}\Big)
					\left(\sum_i \bignorm{\sqrt{\delta_i} a(u_I-u_N)'}_{I_i}^2\right)^{\!1/2}.
		\end{aligned}
		\end{gather}
		Furthermore, we have on the one hand
		\begin{gather*}
			\left|\delta_i \int_{I_i} \eps (u_I-u)'' a (u_I-u_N)' \, dx\right|
				\leq \bignorm{\sqrt{\delta_i}\,\eps(u_I-u)''}_{I_i} \bignorm{\sqrt{\delta_i}a(u_I-u_N)'}_{I_i}
		\end{gather*}
		and on the other hand by the properties of $a$
		\begin{align*}
			\left|\delta_i \int_{I_i} \eps (u_I-u)'' a (u_I-u_N)' \, dx\right|
				& \leq \bignorm{\delta_i \sqrt{\eps} \, a (u_I-u)''}_{I_i} \sqrt{\eps}\snorm{u_I-u_N}_{1,I_i} \\
				& \leq C \norm{\delta_i \sqrt{\eps} \, x(u_I-u)''}_{I_i} \sqrt{\eps} \snorm{u_I-u_N}_{1,I_i}.
		\end{align*}
		Thus, for the first term in~\eqref{ieq:Bstab} the Cauchy Schwarz inequality yields
		\begin{align*}
			& \left|\sum_i \delta_i \int_{I_i} \eps (u_I-u)'' a (u_I-u_N)' \, dx \right|\\
			& \qquad \leq C \left(\sum_i
					\min\left\{ \bignorm{\sqrt{\delta_i}\, \eps (u_I-u)''}_{I_i}^2,
						\norm{\delta_i \sqrt{\eps} \, x(u_I-u)''}_{I_i}^2\right\}\right)^{\!1/2}
					\tnorm{u_I-u_N}_{\mathrm{SD}}.
		\end{align*}
		Combining the above estimates we are done.
	\end{proof}
	
	In Section~\ref{sec:energySunStynes} several interpolation error terms
	have been already studied and bounded. It remains to estimate the last term in Lemma~\ref{le:SDFEMuIuN_Pk}.
	\begin{lemma}
		Let $u$ be the solution of problem~\eqref{prob:intLayer} and $u_I \in V^N$ be its interpolant
		on the piecewise equidistant mesh given by~\eqref{sigma} -- \eqref{eq:hiEll}.
		Suppose that
		\begin{gather}
			\label{ass:delta2}
			0 \leq \delta_i \leq C \min\!\left\{1\, ,\, h_i^2/\eps \right\}\!.
		\end{gather}
		Then
		\begin{gather}
			\label{ieq:secDerIntSunStynes}
			\left(\sum_{i=-N+1}^N \min\left\{ \bignorm{\sqrt{\delta_i}\, \eps (u_I-u)''}_{I_i}^2,
						\norm{\delta_i \sqrt{\eps} \, x(u_I-u)''}_{I_i}^2\right\} \right)^{\!1/2}
				\leq C \max_i\sqrt{\delta_i} \left((\KAK + 1) N^{-1}\right)^{k}.
		\end{gather}
	\end{lemma}
	\begin{proof}
		The proof is similar to the proof of Lemma~\ref{le:IntErrorSunStynes}
		but differs in some details.
		
		Let $x \in (x_{i-1},x_i)$ for some $i$, where $x_i \in (10^{-\KAK},1]$. Then with $\sqrt{\delta_i\eps} \leq C h_i$
		\begin{align*}
			\delta_i \sqrt{\eps} \left| x(u_I-u)''(x) \right|
				& \leq C \big(\sqrt{\delta_i} \, h_i (x_{i-1}+h_i)\big) h_i^{k-1}
					\left(1 + \left( x_{i-1} + \eps^{1/2} \right)^{\lambda - (k+1)}\right) \\
				& \leq C \sqrt{\delta_i} \left((\KAK + 1) N^{-1}\right)^k,
		\end{align*}
		where we used~\eqref{standInt}, \eqref{ieq:innerLayerBounds}, and Lemma~\ref{le:sunStynesh^kBounds}.
		Hence, for $N/(\KAK+1)+1 \leq i \leq N$
		\begin{gather*}
			\norm{\delta_i \sqrt{\eps} \, x(u_I-u)''}_{I_i}^2
				= \int_{I_i} \big( \delta_i \sqrt{\eps} \, x(u_I-u)''(x)\big)^2 dx
				\leq C h_i \delta_i \left((\KAK + 1) N^{-1}\right)^{2k}.
		\end{gather*}
		
		Now, let $x \in (x_{i-1},x_i)$, where $x_i\in(0,10^{-\KAK}]$.
		We consider two different cases for $\sigma$.
		
		First, if $\sigma = \eps^{(1-\lambda/(k+1))/2}$ then as above
		Lemma~\ref{le:sunStynesh^kBounds} yields
		\begin{gather*}
			\delta_i \sqrt{\eps} \left| x(u_I-u)''(x) \right|
				\leq C \sqrt{\delta_i} \left((\KAK + 1) N^{-1}\right)^k
		\end{gather*}
		and therefore
		\begin{gather*}
			\norm{\delta_i \sqrt{\eps} \, x(u_I-u)''}_{I_i}^2
				\leq C h_i \delta_i \left((\KAK + 1) N^{-1}\right)^{2k}
		\end{gather*}
		for $1 \leq i \leq N/(\KAK+1)$.
		
		If $\sigma = N^{-(2k+1)}$ the integral over $(0,x_1)$ can be estimated directly. We have
		\begin{align*}
			\int_0^{x_1} \left(\sqrt{\delta_1}\, \eps (u_I-u)''(x)\right)^2 dx
				& \leq \delta_1 \eps^2 \norm{(u_I-u)''}_{\infty,(0,x_1)}^2 \int_0^{x_1} 1\, dx \\
				& \leq C \delta_1 x_1 \eps^2 \norm{u''}_{\infty,(0,x_1)}^2 \\
				& \leq C \delta_1 x_1 \leq C \delta_1 (\KAK + 1) N^{-2(k+1)}
		\end{align*}
		by~\eqref{standIntInfty}, \eqref{ieq:innerLayerBounds}, and Lemma~\ref{le:sunStynesh^kBounds},
		especially~\eqref{ieq:sunStynesX_1}. For $x \in (x_{i-1},x_i) \subseteq (x_1,10^{-\KAK}]$ we
		use~\eqref{ieq:sunStynes_im1} to obtain
		\begin{align*}
			\delta_i \sqrt{\eps} \left| x(u_I-u)''(x) \right|
				& \leq C \big(\sqrt{\delta_i} \, h_i (x_{i-1}+h_i)\big) h_i^{k-1}
					\left(1 + \left( x_{i-1} + \eps^{1/2} \right)^{\lambda - (k+1)}\right) \\
				& \leq C \sqrt{\delta_i} \left(10^{-\KAK} h_i^k + \left(i-1\right)^{-k} \right),
					\quad \text{for} \quad i=2, \ldots,N/(\KAK+1).
		\end{align*}
		Thus, the inequality~\eqref{ieq:sunStynesX_1} yields for $2\leq i \leq N/(\KAK+1)$
		\begin{align*}
			\norm{\delta_i \sqrt{\eps} \, x(u_I-u)''}_{I_i}^2
				& \leq C \delta_i h_i \left(10^{-\KAK}h_i^k + \left(i-1\right)^{-k} \right)^2 \\
				& \leq C \delta_i x_1 \left(h_i + \left(i-1\right)^{-2} \right) \\
				& \leq C \delta_i \left(h_i + \left(i-1\right)^{-2} \right) (\KAK+1) N^{-2(k+1)}.
		\end{align*}
		
		Summing up the above estimates gives
		\begin{align*}
			& \sum_{i=1}^N \min\left\{ \bignorm{\sqrt{\delta_i}\, \eps (u_I-u)''}_{I_i}^2,
						\norm{\delta_i \sqrt{\eps} \, x(u_I-u)''}_{I_i}^2\right\} \\
			& \qquad \leq C \Bigg( \delta_1 + \!\sum_{i=2}^{N/(\KAK+1)}\delta_i\left(h_i + \left(i-1\right)^{-2}\right) \!\Bigg)
					(\KAK+1) N^{-2(k+1)}
			+ C \sum_{i=1}^N \delta_i h_i \left((\KAK + 1) N^{-1}\right)^{2k} \\
			& \qquad \leq C \max_i \delta_i \left((\KAK + 1) N^{-1}\right)^{2k} \left(1 + N^{-2}\left(1+10^{-\KAK}+\tfrac{\pi^2}{6}\right)\right) \\
			& \qquad \leq C \max_i \delta_i \left((\KAK + 1) N^{-1}\right)^{2k}.
		\end{align*}
		Thanks to symmetry the sum for $i=-N+1,\ldots,0$ can be bounded analogously
		and the proof is completed.
	\end{proof}
	
	The previous estimates enable us to prove an error estimate in the SDFEM-norm.
	\begin{theorem}
		Let $u$ be the solution of~\eqref{prob:intLayer} and $u_N$ the solution
		of~\eqref{dprob:SDFEM} on the piecewise equidistant mesh given
		by~\eqref{sigma} -- \eqref{eq:hiEll}. Furthermore, choose $\delta_i$ such that~\eqref{ass:delta}
		and~\eqref{ass:delta2} are satisfied. Then we have
		\begin{gather*}
			\tnorm{u-u_N}_{\mathrm{SD}} \leq C \left((\KAK + 1) N^{-1}\right)^k \leq C \left(N^{-1}\ln N\right)^k.
		\end{gather*}
	\end{theorem}
	\begin{proof}
		The fact that $a(x)=-xb(x)$ and $\delta_i \leq C$ imply for all $v \in H^1((-1,1))$
		\begin{gather*}
			\left(\sum_{i=-N+1}^N \bignorm{\sqrt{\delta_i}av'}_{I_i}^2\right)^{\!1/2}
				\leq C \max_i \sqrt{\delta_i} \norm{x v'}.
		\end{gather*}
		Using this, the wanted estimate follows easily from the triangle inequality,
		Lemma~\ref{le:SDFEMuIuN_Pk}, \eqref{ieq:energyIntSunStynes}, and~\eqref{ieq:secDerIntSunStynes}.
		The second inequality is an immediately consequence of~\eqref{ieq:KAKtoLN}.
	\end{proof}
	
	\subsection{Some improvements for linear elements}
	\label{sec:someImprove}
	Inspecting the proofs of the last section we see that
	\begin{gather*}
		\left|B_{\mathrm{Stab}}(u_I-u,u_I-u_N)\right|
			\leq C \max_i \sqrt{\delta_i} \left((\KAK+1) N^{-1}\right)^k \tnorm{u_I-u_N}_{\mathrm{SD}}\!.
	\end{gather*}
	For linear elements the intermediate estimate~\eqref{ieq:BstabACterms} can be even improved. Indeed,
	integrating by parts we get
	\begin{align*}
		\int_{I_i} a (u_I-u)' a (u_I-u_N)' \, dx
			& = -\int_{I_i} a^2 (u_I-u)\underbrace{(u_I-u_N)''}_{=0} \, dx
				- \int_{I_i} 2aa'(u_I-u)(u_I-u_N)' \, dx \\
			& \hspace{12em} + \underbrace{\left[a^2(u_I-u)(u_I-u_N)'\right]_{x_{i-1}}^{x_i}}_{=0} \\
			& = -2 \int_{I_i} a'(u_I-u) a(u_I-u_N)' \, dx
	\end{align*}
	and by Cauchy Schwarz' inequality
	\begin{gather*}
		\left|\sum_i \delta_i \int_{I_i} a (u_I-u)' a (u_I-u_N)' \, dx\right|
			\leq C \max_i \sqrt{\delta_i} \norm{u_I-u}
				\left(\sum_i\bignorm{\sqrt{\delta_i} a (u_I-u_N)'}_{I_i}^2\right)^{\!1/2}.
	\end{gather*}
	
	Using the argumentation of~\cite[Theorem~5.1]{SS94} one can prove for linear elements
	\begin{gather*}
		\left|B_\eps(u_I-u,u_I-u_N)\right|
			\leq C \left((\KAK+1)N^{-1}\right)^{3/2} \norm{u_I-u_N}
			\leq C \left(N^{-1}\ln N\right)^{3/2} \norm{u_I-u_N}.
	\end{gather*}
	Note that the occurring logarithmic factors originate from an estimate
	of the form $(\KAK+1)\leq C \ln N$.
	
	In summary we get the following result.
	\begin{theorem}
		\label{th:sdfemSunStynes}
		Let $u$ be the solution of~\eqref{prob:intLayer}, $u_N \in V^N$ $(k = 1)$ the solution
		of~\eqref{dprob:SDFEM}, and $u_I \in V^N$ $(k=1)$ the interpolant of $u$ on the piecewise
		equidistant mesh given by~\eqref{sigma} -- \eqref{eq:hiEll}. Furthermore, choose
		$\delta_i$ such that
		\begin{gather*}
			0 \leq \delta_i \leq \min\left\{ \frac{\gamma}{2 \norm{c}_{\infty}^2}, \frac{h_i^2}{\eps}, (\KAK+1)N^{-1} \right\}.
		\end{gather*}
		Then we have
		\begin{gather*}
			\norm{u-u_N} + \tnorm{u_I-u_N}_{\mathrm{SD}}
				\leq C \left((\KAK + 1) N^{-1}\right)^{3/2} \leq C \left(N^{-1}\ln N\right)^{3/2}.
		\end{gather*}
	\end{theorem}
	\begin{proof}
		Revise the proof of Lemma~\ref{le:SDFEMuIuN_Pk} using the improved estimates of this section. Then
		invoke~\eqref{ieq:secDerIntSunStynes} and~\eqref{ieq:L2IntSunStynes} to complete the proof.
	\end{proof}
	
	\begin{remark}
		For linear Galerkin FEM in~\cite[Theorem~5.1]{SS94} the $L^2$-norm
		estimate
		\begin{gather*}
			\norm{u-u_N} \leq C\left(N^{-1}\ln N\right)^{3/2}
		\end{gather*}
		is already proven for the piecewise equidistant mesh.
	\end{remark}
	
	\section{Numerical experiments}
	\label{sec:numExp}
	In this section we present some numerical results to illustrate the theoretical findings
	of the previous sections. As in~\cite{Bec16ArXivLis} we study a test problem taken
	from~\cite{SS94} whose solution exhibits typical interior ``cusp''-type layer behaviour.
	
	\begin{example}[see~\cite{SS94}]
		\label{test:SunStynes}
		We consider the singularly perturbed turning point problem
		\begin{align*}
			- \eps u'' - x(1+x^2) u' + \lambda(1+x^3)u & = f, \qquad \text{for} \quad x \in (-1,1), \\
			u(-1)=u(1) & = 0,
		\end{align*}
		where the right-hand side $f(x)$ is chosen such that the solution $u(x)$ is given by
		\begin{gather*}
			u(x) = \left(x^2+\eps\right)^{\lambda/2} + x\left(x^2+\eps\right)^{(\lambda-1)/2}
				- \left(1+\eps\right)^{\lambda/2}\left(1+x\left(1+\eps\right)^{-1/2}\right).
		\end{gather*}
		Note that the parameter $\lambda$ in the problem coincides with the quantity
		$\bar{\lambda}=c(0)/|a'(0)|$.
	\end{example}
	
	All computations were performed using a FEM-code based on $\mathbb{SOFE}$ by Lars
	Ludwig~\cite{SOFE}. Motivated by our error estimates we calculate the
	convergence rates by
	\begin{align*}
		r   & = \left(\ln E_{\eps,N}-\ln E_{\eps,2N}\right)/\ln 2
	\end{align*}
	for given errors $E_{\eps,N}$. 
	In order to ensure that the used meshes consists of exactly $2N$ mesh intervals,
	which is presumed in this formula, we adjust the mesh like in~\cite[Section~6]{SS94}:
	
	Let $N_0 = N-(\KAK+1)n_0$ where $n_0=\lfloor N/(\KAK+1)\rfloor$. Then we uniformly
	partition each of the subintervals $(0,10^{-\KAK}], \ldots, (10^{-N_0-1},10^{-N_0}]$
	by $n_0$ points and each of the remaining subintervals
	$(10^{-N_0},10^{-N_0+1}], \ldots, \allowbreak (10^{-1},1]$ by $n_0+1$ points. The arising mesh
	is still piecewise equidistant and has exactly $2N$ mesh intervals.
	
	For the streamline-diffusion finite element method we choose the stabilisation
	parameter as
	\begin{gather*}
		\delta_i = C_0 \min \big\{h_i^2/\eps, h_i\big\}
	\end{gather*}
	which is the standard choice, see e.g.~\cite[p.~87]{RST08}. Although it is not
	necessary to have $\delta_i \leq C_0 h_i$ by theory, numerical tests suggest
	to favour this definition. We use $C_0=1$ for computations.
	
\begin{figure}
		\begin{center}
			\input{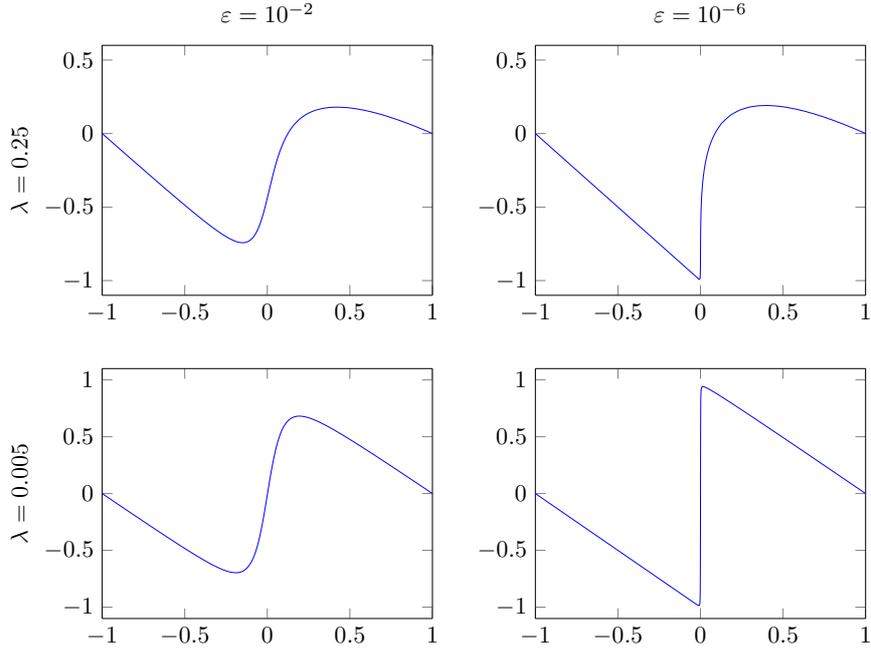}
			\caption{Numerical solutions of $P_2$-FEM on a mesh with $N=128$ applied to
				Example~\ref{test:SunStynes} with $\eps \in \{10^{-2},10^{-6}\}$ (left/right)
				and $\lambda \in \{0.25, 0.005\}$ (top/bottom row).}
			\label{fig:numsol}
		\end{center}
	\end{figure}

	Numerical solutions of Example~\ref{test:SunStynes} are displayed in
	Figure~\ref{fig:numsol} for various values of $\eps$ and $\lambda$.
	Here $P_2$-FEM was applied on a mesh with $N=128$.
	
\begin{figure}
		\begin{center}
%
%
%
\begin{tikzpicture}[scale=0.9]

\begin{axis}[%
width=4.4in,
height=3.5in,
at={(0.809in,0.513in)},
scale only axis,
separate axis lines,
every outer x axis line/.append style={black},
every x tick label/.append style={font=\color{black}},
xmode=log,
xmin=10,
xmax=10000,
xtick={  32,   64,  128,  256,  512, 1024, 2048, 4096},
xticklabels={  32,   64,  128,  256,  512, 1024, 2048, 4096},
xminorticks=true,
xlabel={number of mesh intervals $\KAO(N)$},
every outer y axis line/.append style={black},
every y tick label/.append style={font=\color{black}},
ymode=log,
ymin=1e-15,
ymax=1,
yminorticks=true,
ylabel={energy norm error/SDFEM-norm error},
axis background/.style={fill=white},
legend style={at={(0.03,0.03)},anchor=south west,legend cell align=left,align=left,draw=black}
]
\addplot [color=blue,dashed,mark=triangle*,mark options={solid}]
  table[row sep=crcr]{%
32	0.000696355758031245\\
64	0.000339146857716946\\
128	0.000163706521233724\\
256	8.04754284204038e-05\\
512	3.97229372277397e-05\\
1024	1.98313102234079e-05\\
2048	9.88577167729452e-06\\
4096	4.94233941153684e-06\\
};
\addlegendentry{$P_1$-FEM};
\addplot [color=blue,dashed,mark=square,mark options={solid}]
  table[row sep=crcr]{%
32	0.000860151241151813\\
64	0.000406575126864472\\
128	0.000180228975678878\\
256	8.55624943568482e-05\\
512	4.09892821292261e-05\\
1024	2.01650012167116e-05\\
2048	9.96898880591745e-06\\
4096	4.96344411075546e-06\\
};
\addlegendentry{$P_1$-SDFEM};

\addplot [color=blue,solid,mark=triangle*,mark options={solid}]
  table[row sep=crcr]{%
32	0.000223620532498402\\
64	7.04168966100957e-05\\
128	1.7480053133063e-05\\
256	4.45529100262357e-06\\
512	1.09790358486655e-06\\
1024	2.73300673861802e-07\\
2048	6.72887123201344e-08\\
4096	1.66990078016726e-08\\
};
\addlegendentry{$P_2$-FEM};
\addplot [color=blue,solid,mark=square,mark options={solid}]
  table[row sep=crcr]{%
32	0.000212856400063806\\
64	6.37458332033132e-05\\
128	1.48520037469274e-05\\
256	3.62240803624461e-06\\
512	8.64027845018387e-07\\
1024	2.12581122594421e-07\\
2048	5.23191976095532e-08\\
4096	1.30055434018444e-08\\
};
\addlegendentry{$P_2$-SDFEM};

\addplot [color=blue,dotted,mark=triangle*,mark options={solid}]
  table[row sep=crcr]{%
32	4.49057950325087e-05\\
64	6.94582501693579e-06\\
128	7.64539591245552e-07\\
256	9.04583193550302e-08\\
512	1.08292063480329e-08\\
1024	1.3438111977864e-09\\
2048	1.6638789211599e-10\\
4096	2.07852550129963e-11\\
};
\addlegendentry{$P_3$-FEM};
\addplot [color=blue,dotted,mark=square,mark options={solid}]
  table[row sep=crcr]{%
32	4.77091134005691e-05\\
64	7.74163847177039e-06\\
128	8.3026049111342e-07\\
256	9.79162595985581e-08\\
512	1.13192779651683e-08\\
1024	1.37945076858667e-09\\
2048	1.6861951585935e-10\\
4096	2.09298865978016e-11\\
};
\addlegendentry{$P_3$-SDFEM};

\addplot [color=blue,dashdotted,mark=triangle*,mark options={solid}]
  table[row sep=crcr]{%
32	9.0802243142672e-06\\
64	7.05835199072986e-07\\
128	4.37458213747794e-08\\
256	3.1356912683582e-09\\
512	1.96626947485262e-10\\
1024	1.22507584597869e-11\\
2048	7.53705582607914e-13\\
4096	7.02978533850232e-14\\
};
\addlegendentry{$P_4$-FEM};
\addplot [color=blue,dashdotted,mark=square,mark options={solid}]
  table[row sep=crcr]{%
32	9.24483613812883e-06\\
64	7.18526974189683e-07\\
128	3.90065454165416e-08\\
256	2.42826791691573e-09\\
512	1.35470756420519e-10\\
1024	7.8197240150258e-12\\
2048	4.57782718149712e-13\\
4096	6.27600451968103e-14\\
};
\addlegendentry{$P_4$-SDFEM};

\addplot [color=blue,dashed,forget plot]
  table[row sep=crcr]{%
32	0.000375\\
64	0.0001875\\
128	9.375e-05\\
256	4.6875e-05\\
512	2.34375e-05\\
1024	1.171875e-05\\
2048	5.859375e-06\\
4096	2.9296875e-06\\
};
\addplot [color=blue,solid,forget plot]
  table[row sep=crcr]{%
32	0.0001171875\\
64	2.9296875e-05\\
128	7.32421875e-06\\
256	1.8310546875e-06\\
512	4.57763671875e-07\\
1024	1.1444091796875e-07\\
2048	2.86102294921875e-08\\
4096	7.15255737304687e-09\\
};
\addplot [color=blue,dotted,forget plot]
  table[row sep=crcr]{%
32	2.13623046875e-05\\
64	2.6702880859375e-06\\
128	3.33786010742187e-07\\
256	4.17232513427734e-08\\
512	5.21540641784668e-09\\
1024	6.51925802230835e-10\\
2048	8.14907252788544e-11\\
4096	1.01863406598568e-11\\
};
\addplot [color=blue,dashdotted,forget plot]
  table[row sep=crcr]{%
32	4.76837158203125e-06\\
64	2.98023223876953e-07\\
128	1.86264514923096e-08\\
256	1.16415321826935e-09\\
512	7.27595761418343e-11\\
1024	4.54747350886464e-12\\
2048	2.8421709430404e-13\\
4096	1.77635683940025e-14\\
};
\end{axis}
\end{tikzpicture}%
			\caption{Energy norm error for $P_k$-FEM and SDFEM-norm error for $P_k$-SDFEM,
			$k = 1, \ldots, 4$, applied to Example~\ref{test:SunStynes} with $\eps = 10^{-10}$
			and $\lambda = 0.005$. Reference curves of the form $\KAO(N^{-k})$.}
			\label{fig:SdFem}
		\end{center}
	\end{figure}
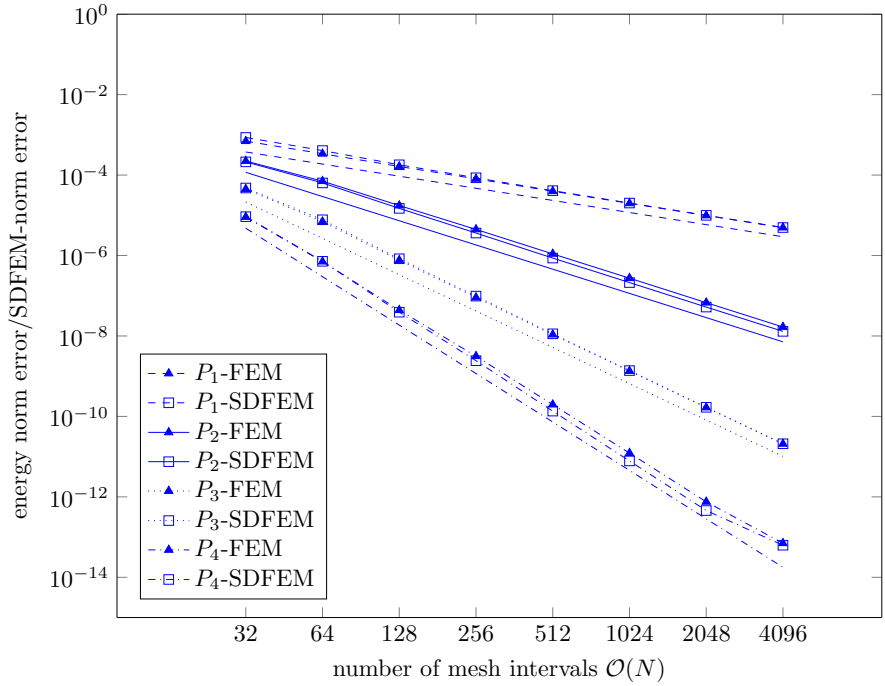

	We plot the energy norm error for $P_k$-FEM and the SDFEM-norm error
	for $P_k$-SDFEM, $k = 1, \ldots, 4$, in Figure~\ref{fig:SdFem} where the
	methods were applied to Example~\ref{test:SunStynes} with $\eps = 10^{-10}$
	and $\lambda = 0.005$. The expected orders of convergence, cf.
	Theorem~\ref{th:energyErrorSunStynes}, Theorem~\ref{th:sdfemSunStynes},
	can be clearly seen. The magnitude of the errors is similar for both
	methods. Recall that the SDFEM-norm is stronger than the energy
	norm. The numerical results also suggest that the errors are uniform
	with respect to $\eps$. They stay stable for small $\eps$, see
	Table~\ref{tab:energy} for FEM and Table~\ref{tab:sdfem} for SDFEM.
	
	\begin{table}
	\begin{center}
	\begin{tabular}{ l *{8}{c} } \toprule
			& \multicolumn{2}{c}{$P_1$-elements} & \multicolumn{2}{c}{$P_2$-elements} &
				\multicolumn{2}{c}{$P_3$-elements} & \multicolumn{2}{c}{$P_4$-elements} \\
		\cmidrule(l{0.7ex}r{0.7ex}){2-3} \cmidrule(l{0.7ex}r{0.7ex}){4-5}
		\cmidrule(l{0.7ex}r{0.7ex}){6-7} \cmidrule(l{0.7ex}r{0.7ex}){8-9}
		$\eps$ \textbf{\textbackslash} $N$ & 512 & 1024 & 512 & 1024 & 512 & 1024 & 512 & 1024 \\ 
		\cmidrule(l{0.7ex}r{0.7ex}){1-1} \cmidrule(l{0.7ex}r{0.7ex}){2-2} \cmidrule(l{0.7ex}r{0.7ex}){3-3}
		\cmidrule(l{0.7ex}r{0.7ex}){4-4} \cmidrule(l{0.7ex}r{0.7ex}){5-5} \cmidrule(l{0.7ex}r{0.7ex}){6-6}
		\cmidrule(l{0.7ex}r{0.7ex}){7-7} \cmidrule(l{0.7ex}r{0.7ex}){8-8} \cmidrule(l{0.7ex}r{0.7ex}){9-9}
		$1$ 	   & 9.71e-04 & 4.85e-04 & 7.49e-07 & 1.87e-07 & 9.17e-10 & 1.15e-10 & 6.53e-13 & 2.10e-12 \\ 
		$10^{- 2}$ & 1.38e-03 & 6.89e-04 & 9.60e-06 & 2.40e-06 & 4.10e-08 & 5.12e-09 & 1.34e-10 & 8.39e-12 \\ 
		$10^{- 4}$ & 6.45e-04 & 3.24e-04 & 6.73e-06 & 1.69e-06 & 4.32e-08 & 5.44e-09 & 2.20e-10 & 1.39e-11 \\ 
		$10^{- 6}$ & 2.69e-04 & 1.35e-04 & 3.76e-06 & 9.40e-07 & 3.22e-08 & 4.02e-09 & 2.20e-10 & 1.38e-11 \\ 
		$10^{- 8}$ & 1.06e-04 & 5.26e-05 & 1.90e-06 & 4.70e-07 & 1.99e-08 & 2.45e-09 & 1.82e-10 & 1.11e-11 \\ 
		$10^{-10}$ & 3.97e-05 & 1.98e-05 & 1.10e-06 & 2.73e-07 & 1.08e-08 & 1.34e-09 & 1.97e-10 & 1.23e-11 \\ 
		$10^{-12}$ & 1.47e-05 & 7.25e-06 & 1.06e-06 & 2.65e-07 & 5.51e-09 & 6.69e-10 & 3.21e-10 & 2.02e-11 \\ 
		$10^{-14}$ & 7.48e-06 & 2.91e-06 & 1.33e-06 & 3.34e-07 & 5.05e-09 & 4.17e-10 & 5.59e-10 & 3.51e-11 \\ 
		\bottomrule
	\end{tabular}
	\caption{Energy norm error for $P_k$-FEM, $k = 1, \ldots, 4$, applied
	to Example~\ref{test:SunStynes} with certain $\eps$ and $\lambda = 0.005$.}
	\label{tab:energy}
	\end{center}
	\end{table}
	
	\begin{table}
	\begin{center}
	\begin{tabular}{ l *{8}{c} } \toprule
			& \multicolumn{2}{c}{$P_1$-elements} & \multicolumn{2}{c}{$P_2$-elements} &
				\multicolumn{2}{c}{$P_3$-elements} & \multicolumn{2}{c}{$P_4$-elements} \\
		\cmidrule(l{0.7ex}r{0.7ex}){2-3} \cmidrule(l{0.7ex}r{0.7ex}){4-5}
		\cmidrule(l{0.7ex}r{0.7ex}){6-7} \cmidrule(l{0.7ex}r{0.7ex}){8-9}
		$\eps$ \textbf{\textbackslash} $N$ & 512 & 1024 & 512 & 1024 & 512 & 1024 & 512 & 1024 \\ 
		\cmidrule(l{0.7ex}r{0.7ex}){1-1} \cmidrule(l{0.7ex}r{0.7ex}){2-2} \cmidrule(l{0.7ex}r{0.7ex}){3-3}
		\cmidrule(l{0.7ex}r{0.7ex}){4-4} \cmidrule(l{0.7ex}r{0.7ex}){5-5} \cmidrule(l{0.7ex}r{0.7ex}){6-6}
		\cmidrule(l{0.7ex}r{0.7ex}){7-7} \cmidrule(l{0.7ex}r{0.7ex}){8-8} \cmidrule(l{0.7ex}r{0.7ex}){9-9}
		$1$		   & 9.71e-04 & 4.85e-04 & 7.49e-07 & 1.87e-07 & 9.17e-10 & 1.15e-10 & 7.11e-13 & 2.43e-12 \\ 
		$10^{- 2}$ & 1.39e-03 & 6.90e-04 & 1.02e-05 & 2.44e-06 & 4.98e-08 & 5.39e-09 & 2.26e-10 & 9.54e-12 \\ 
		$10^{- 4}$ & 6.46e-04 & 3.24e-04 & 6.76e-06 & 1.70e-06 & 4.40e-08 & 5.56e-09 & 2.68e-10 & 1.90e-11 \\ 
		$10^{- 6}$ & 2.69e-04 & 1.35e-04 & 3.76e-06 & 9.39e-07 & 3.22e-08 & 4.02e-09 & 2.19e-10 & 1.37e-11 \\ 
		$10^{- 8}$ & 1.06e-04 & 5.27e-05 & 1.85e-06 & 4.58e-07 & 2.00e-08 & 2.46e-09 & 1.73e-10 & 1.05e-11 \\ 
		$10^{-10}$ & 4.10e-05 & 2.02e-05 & 8.64e-07 & 2.13e-07 & 1.13e-08 & 1.38e-09 & 1.35e-10 & 7.82e-12 \\ 
		$10^{-12}$ & 1.95e-05 & 8.58e-06 & 4.89e-07 & 1.09e-07 & 7.98e-09 & 8.53e-10 & 1.70e-10 & 8.37e-12 \\ 
		$10^{-14}$ & 1.72e-05 & 6.28e-06 & 5.01e-07 & 9.67e-08 & 1.02e-08 & 9.24e-10 & 3.02e-10 & 1.38e-11 \\ 
		\bottomrule
	\end{tabular}
	\caption{SDFEM-norm error for $P_k$-SDFEM, $k = 1, \ldots, 4$, applied
	to Example~\ref{test:SunStynes} with certain $\eps$ and $\lambda = 0.005$.}
	\label{tab:sdfem}
	\end{center}
	\end{table}
	
	In Figure~\ref{fig:lambda025} and Figure~\ref{fig:lambda0005} the errors
	of both methods are plotted for $P_2$-elements applied to Example~\ref{test:SunStynes}
	with $\eps = 10^{-6}$, $\lambda=0.25$ and $\lambda = 0.005$,
	respectively, on a mesh with $N=128$. While for the finite element method
	the error is clearly oscillating in $(0,1)$, this behaviour is damped and
	even prevented in a wide range of the interval when the stabilisation
	technique is used.
	
	{
	\begin{figure}
		\begin{center}
			\input{lambda025N128fem}
			\input{lambda025N128sdfem}
			\caption{Error for $P_2$-FEM (left) and $P_2$-SDFEM (right) applied to
			Example~\ref{test:SunStynes} with $\eps = 10^{-6}$ and $\lambda = 0.25$
			on a mesh with $N = 128$.}
			\label{fig:lambda025}
		\end{center}
	\end{figure}
	\begin{figure}
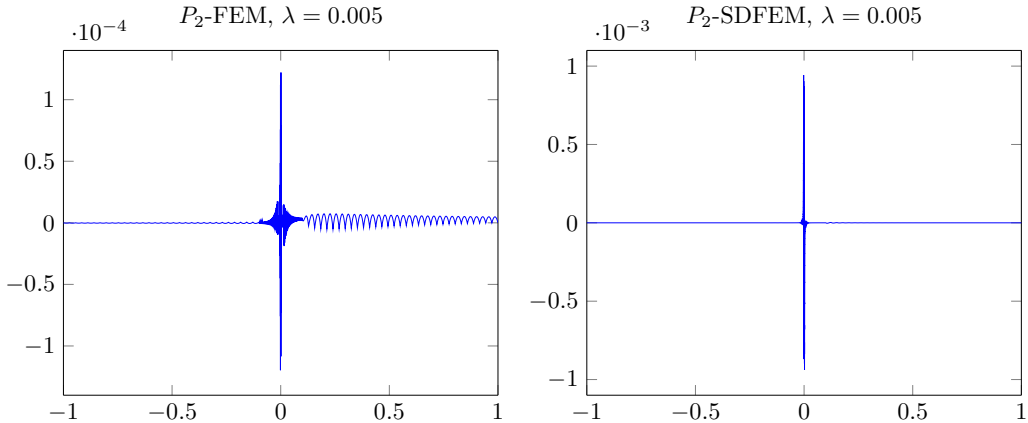

		\begin{center}
			\input{lambda0005N128fem}
			\input{lambda0005N128sdfem}
			\caption{Error for $P_2$-FEM (left) and $P_2$-SDFEM (right) applied to
			Example~\ref{test:SunStynes} with $\eps = 10^{-6}$ and $\lambda = 0.005$
			on a mesh with $N = 128$.}
			\label{fig:lambda0005}
		\end{center}
	\end{figure}
	}
	
	For linear SDFEM the errors in SDFEM-norm, energy norm, and $L^2$-norm can be
	compared in Table~\ref{tab:l2norm}. Here the method was applied to Example~\ref{test:SunStynes}
	with parameters $\eps = 10^{-10}$ and $\lambda = 0.005$. The calculated convergence
	rates coincide with the rates expected from standard theory for non singularly
	perturbed problems. So we have first order convergence in the SDFEM-norm and second
	order convergence in the $L^2$-norm. But recall that in Section~\ref{sec:someImprove}
	we were able to prove the convergence rate $3/2$ only.
	
	\begin{table}
	\begin{center}
	\begin{tabular}{ r *{6}{c} } \toprule
			& \multicolumn{2}{c}{$\tnorm{u-u_N}_{\mathrm{SD}}$} & \multicolumn{2}{c}{$\tnorm{u-u_N}_\eps$}
			& \multicolumn{2}{c}{$\norm{u-u_N}$} \\
		$N$ & error & rates & error & rates & error & rates \\ 
		\cmidrule(l{0.7ex}r{0.7ex}){1-1} \cmidrule(l{0.7ex}r{0.7ex}){2-3}
		\cmidrule(l{0.7ex}r{0.7ex}){4-5} \cmidrule(l{0.7ex}r{0.7ex}){6-7} 
		   8 & 2.23e-02 & 2.584    & 2.16e-02 & 2.601    & 2.16e-02 & 3.177    \\
		  16 & 3.71e-03 & 2.110    & 3.57e-03 & 2.493    & 2.39e-03 & 2.989    \\ 
		  32 & 8.60e-04 & 1.081    & 6.33e-04 & 0.943    & 3.01e-04 & 1.699    \\ 
		  64 & 4.07e-04 & 1.174    & 3.30e-04 & 1.045    & 9.27e-05 & 1.982    \\ 
		 128 & 1.80e-04 & 1.075    & 1.60e-04 & 0.995    & 2.35e-05 & 2.006    \\ 
		 256 & 8.56e-05 & 1.062    & 8.02e-05 & 1.016    & 5.84e-06 & 2.086    \\ 
		 512 & 4.10e-05 & 1.023    & 3.96e-05 & 1.000    & 1.38e-06 & 2.059    \\ 
		1024 & 2.02e-05 & 1.016    & 1.98e-05 & 1.004    & 3.30e-07 & 2.055    \\ 
		2048 & 9.97e-06 & 1.006    & 9.88e-06 & 1.000    & 7.95e-08 & 2.032    \\ 
		\bottomrule
	\end{tabular}
	\caption{SDFEM-norm, energy norm, and $L^2$-norm error for linear SDFEM applied
	to Example~\ref{test:SunStynes} with $\eps = 10^{-10}$ and $\lambda = 0.005$.}
	\label{tab:l2norm}
	\end{center}
	\end{table}
	
	Furthermore, we want to present some numerical computations for $\lambda = 0.25$.
	In Table~\ref{tab:energyL2norm1_2} and Table~\ref{tab:energyL2norm3_4} the energy
	and $L^2$-norm errors are given together with the associated convergence rates
	for certain $\eps$ and $P_k$-FEM, $k=1,\ldots,4$. In the studied range of $N$
	the results qualitatively differ from the results for smaller $\lambda$
	(before $\lambda = 0.005$ was studied). On the one hand for $\eps \geq 10^{-4}$
	the expected convergence behaviour can be seen. Otherwise for smaller $\eps$
	the $L^2$-norm part dominates the energy norm error which may suggest that
	this norm is too weak. Here we also have to differentiate between odd and
	even element orders. For $k=1,3$ we obtain the convergence order $k+1$ in
	the $L^2$-norm.	So for small $\eps$ also the rate in the $\tnorm{\cdot}_\eps$-norm
	is calculated to be $k+1$ which surpasses the usual expectations. Otherwise
	for $k=2,4$ we see the expected rate $k$ in the energy norm and when
	$\eps \leq 10^{-8}$ the same rate also for $\norm{u-u_N}$ where one
	would rather expect $k+1$ as obtained for larger $\eps$.
	
	Finally, we want to point out that the detailed structure of the error estimate becomes
	visible for $\lambda=0.25$ at least for $P_2$-elements. In order to check this we
	additionally calculate the ratio of the numerically computed error to the proven
	error bound which is given by $\tnorm{u-u_N}_\eps \cdot 100 {\big(N/(\KAK+1)\big)^2}$,
	see Table~\ref{tab:KAKfem2}. Excluding the results for $\eps \geq 10^{-2}$
	the ratio is nearly independent of $\eps$. Especially in the lower right corner
	(for $N \geq 256$ and $\eps \leq 10^{-10}$) the ratio is even nearly constant.
	
	\begin{table}
	\begin{center}
	\begin{tabular}{ l *{8}{c} } \toprule
			& \multicolumn{4}{c}{$P_1$-elements} & \multicolumn{4}{c}{$P_2$-elements} \\
		\cmidrule(l{0.7ex}r{0.7ex}){2-5} \cmidrule(l{0.7ex}r{0.7ex}){6-9}
			& \multicolumn{2}{c}{$\tnorm{u-u_N}_\eps$} & \multicolumn{2}{c}{$\norm{u-u_N}$}
			& \multicolumn{2}{c}{$\tnorm{u-u_N}_\eps$} & \multicolumn{2}{c}{$\norm{u-u_N}$} \\[0.5ex]
		$\eps$ \textbf{\textbackslash} $N$ & 512 & 1024 & 512 & 1024 & 512 & 1024 & 512 & 1024 \\
		\cmidrule(l{0.7ex}r{0.7ex}){1-1}
		\cmidrule(l{0.7ex}r{0.7ex}){2-3} \cmidrule(l{0.7ex}r{0.7ex}){4-5}
		\cmidrule(l{0.7ex}r{0.7ex}){6-7} \cmidrule(l{0.7ex}r{0.7ex}){8-9}
		$1$        & 8.06e-04 & 4.03e-04 & 8.65e-07 & 2.16e-07 & 5.94e-07 & 1.49e-07 & 3.22e-10 & 4.03e-11 \\ 
		           & 1.000    &          & 2.000    &          & 2.000    &          & 3.000    &          \\ 
		$10^{- 2}$ & 7.80e-04 & 3.90e-04 & 5.63e-06 & 1.41e-06 & 4.39e-06 & 1.10e-06 & 2.38e-08 & 2.97e-09 \\ 
		           & 1.000    &          & 2.000    &          & 2.000    &          & 3.000    &          \\ 
		$10^{- 4}$ & 2.18e-04 & 1.09e-04 & 5.39e-06 & 1.30e-06 & 1.93e-06 & 4.70e-07 & 6.27e-08 & 6.74e-09 \\ 
		           & 0.998    &          & 2.052    &          & 2.037    &          & 3.218    &          \\ 
		$10^{- 6}$ & 5.63e-05 & 2.67e-05 & 1.89e-05 & 4.03e-06 & 4.43e-06 & 1.01e-06 & 3.46e-06 & 4.85e-07 \\ 
		           & 1.076    &          & 2.225    &          & 2.139    &          & 2.834    &          \\ 
		$10^{- 8}$ & 3.28e-05 & 9.85e-06 & 3.05e-05 & 7.87e-06 & 7.61e-06 & 1.83e-06 & 7.57e-06 & 1.80e-06 \\ 
		           & 1.735    &          & 1.956    &          & 2.055    &          & 2.069    &          \\ 
		$10^{-10}$ & 4.84e-05 & 1.11e-05 & 4.83e-05 & 1.11e-05 & 1.14e-05 & 2.86e-06 & 1.14e-05 & 2.85e-06 \\ 
		           & 2.121    &          & 2.128    &          & 1.992    &          & 1.992    &          \\ 
		$10^{-12}$ & 5.58e-05 & 1.48e-05 & 5.58e-05 & 1.48e-05 & 1.55e-05 & 3.91e-06 & 1.55e-05 & 3.91e-06 \\ 
		           & 1.915    &          & 1.916    &          & 1.986    &          & 1.986    &          \\ 
		$10^{-14}$ & 8.88e-05 & 2.23e-05 & 8.88e-05 & 2.23e-05 & 2.06e-05 & 5.16e-06 & 2.06e-05 & 5.16e-06 \\ 
		           & 1.992    &          & 1.992    &          & 1.998    &          & 1.998    &          \\ 
		\bottomrule
	\end{tabular}
	\caption{Energy norm and $L^2$-norm error for $P_k$-FEM, $k = 1, 2$, applied
	to Example~\ref{test:SunStynes} with certain $\eps$ and $\lambda = 0.25$.
	Associated convergence rates.}
	\label{tab:energyL2norm1_2}
	\end{center}
	\end{table}
	
	\begin{table}
	\begin{center}
	\begin{tabular}{ l *{8}{c} } \toprule
			& \multicolumn{4}{c}{$P_3$-elements} & \multicolumn{4}{c}{$P_4$-elements} \\
		\cmidrule(l{0.7ex}r{0.7ex}){2-5} \cmidrule(l{0.7ex}r{0.7ex}){6-9}
			& \multicolumn{2}{c}{$\tnorm{u-u_N}_\eps$} & \multicolumn{2}{c}{$\norm{u-u_N}$}
			& \multicolumn{2}{c}{$\tnorm{u-u_N}_\eps$} & \multicolumn{2}{c}{$\norm{u-u_N}$} \\[0.5ex]
		$\eps$ \textbf{\textbackslash} $N$ & 512 & 1024 & 512 & 1024 & 512 & 1024 & 512 & 1024 \\
		\cmidrule(l{0.7ex}r{0.7ex}){1-1}
		\cmidrule(l{0.7ex}r{0.7ex}){2-3} \cmidrule(l{0.7ex}r{0.7ex}){4-5}
		\cmidrule(l{0.7ex}r{0.7ex}){6-7} \cmidrule(l{0.7ex}r{0.7ex}){8-9}
		$1$        & 6.72e-10 & 8.40e-11 & 3.92e-13 & 2.51e-13 & 7.77e-12 & 2.83e-11 & 3.70e-12 & 1.36e-11 \\ 
		           & 3.000    &          & 0.643    &          & -1.868   &          & -1.875   &          \\ 
		$10^{- 2}$ & 1.66e-08 & 2.08e-09 & 6.15e-11 & 4.59e-12 & 6.70e-11 & 4.24e-12 & 1.93e-13 & 3.93e-13 \\ 
		           & 3.000    &          & 3.745    &          & 3.982    &          & -1.028   &          \\ 
		$10^{- 4}$ & 1.21e-08 & 1.49e-09 & 3.24e-10 & 2.17e-11 & 9.80e-11 & 6.02e-12 & 3.08e-12 & 1.17e-13 \\ 
		           & 3.013    &          & 3.901    &          & 4.025    &          & 4.726    &          \\ 
		$10^{- 6}$ & 6.36e-09 & 6.98e-10 & 1.83e-09 & 6.78e-11 & 3.82e-10 & 1.93e-11 & 1.51e-10 & 3.47e-12 \\ 
		           & 3.188    &          & 4.754    &          & 4.309    &          & 5.443    &          \\ 
		$10^{- 8}$ & 7.19e-09 & 5.29e-10 & 6.80e-09 & 4.57e-10 & 9.38e-10 & 5.72e-11 & 9.18e-10 & 5.46e-11 \\ 
		           & 3.764    &          & 3.897    &          & 4.036    &          & 4.073    &          \\ 
		$10^{-10}$ & 1.68e-08 & 8.55e-10 & 1.67e-08 & 8.50e-10 & 2.10e-09 & 1.33e-10 & 2.10e-09 & 1.33e-10 \\ 
		           & 4.293    &          & 4.299    &          & 3.978    &          & 3.979    &          \\ 
		$10^{-12}$ & 2.39e-08 & 1.52e-09 & 2.39e-08 & 1.52e-09 & 3.91e-09 & 2.50e-10 & 3.91e-09 & 2.50e-10 \\ 
		           & 3.973    &          & 3.973    &          & 3.969    &          & 3.969    &          \\ 
		$10^{-14}$ & 5.56e-08 & 3.56e-09 & 5.56e-08 & 3.56e-09 & 6.93e-09 & 4.36e-10 & 6.93e-09 & 4.36e-10 \\ 
		           & 3.967    &          & 3.967    &          & 3.991    &          & 3.991    &          \\ 
		\bottomrule
	\end{tabular}
	\caption{Energy norm and $L^2$-norm error for $P_k$-FEM, $k = 3, 4$, applied
	to Example~\ref{test:SunStynes} with certain $\eps$ and $\lambda = 0.25$.
	Associated convergence rates.}
	\label{tab:energyL2norm3_4}
	\end{center}
	\end{table}
	
	\begin{table}
	\begin{center}
	\begin{tabular}{ l *{10}{c} } \toprule
		$\eps$ \textbf{\textbackslash} $N$ & 8 & 16 & 32 & 64 & 128 & 256 & 512 & 1024 & 2048 & 4096 \\
		\cmidrule(l{0.7ex}r{0.7ex}){1-1} \cmidrule(l{0.7ex}r{0.7ex}){2-2} \cmidrule(l{0.7ex}r{0.7ex}){3-3}
		\cmidrule(l{0.7ex}r{0.7ex}){4-4} \cmidrule(l{0.7ex}r{0.7ex}){5-5} \cmidrule(l{0.7ex}r{0.7ex}){6-6}
		\cmidrule(l{0.7ex}r{0.7ex}){7-7} \cmidrule(l{0.7ex}r{0.7ex}){8-8} \cmidrule(l{0.7ex}r{0.7ex}){9-9}
		\cmidrule(l{0.7ex}r{0.7ex}){10-10} \cmidrule(l{0.7ex}r{0.7ex}){11-11}
		$1$        &  3.82 &  3.88 &  3.89 &  3.89 &  3.90 &  3.90 &  3.90 &  3.90 &  3.90 &  3.90 \\ 
		$10^{- 1}$ &  8.99 & 10.05 & 10.72 & 10.88 & 10.92 & 10.93 & 10.94 & 10.94 & 10.94 & 10.94 \\ 
		$10^{- 2}$ & 20.95 & 26.38 & 28.34 & 28.68 & 28.74 & 28.75 & 28.76 & 28.76 & 28.76 & 28.76 \\ 
		$10^{- 3}$ &  5.58 &  7.49 &  9.35 &  8.66 &  8.47 &  8.28 &  8.31 &  8.30 &  8.31 &  8.31 \\ 
		$10^{- 4}$ &  5.15 &  7.06 &  8.24 &  8.50 &  7.40 &  6.21 &  5.61 &  5.47 &  5.41 &  5.41 \\ 
		$10^{- 5}$ &  5.10 &  6.64 &  7.23 &  7.33 &  7.39 &  7.08 &  5.52 &  3.48 &  2.27 &  1.78 \\ 
		$10^{- 6}$ &  5.11 &  6.70 &  7.64 &  7.78 &  7.46 &  7.33 &  7.26 &  6.59 &  4.60 &  2.69 \\ 
		$10^{- 7}$ &  3.28 &  4.74 &  6.38 &  7.76 &  7.84 &  7.64 &  7.33 &  7.26 &  7.18 &  6.49 \\ 
		$10^{- 8}$ &  3.29 &  4.75 &  6.44 &  7.87 &  7.98 &  8.01 &  7.98 &  7.69 &  7.32 &  7.24 \\ 
		$10^{- 9}$ &  2.34 &  4.85 &  5.95 &  7.58 &  7.83 &  8.22 &  8.17 &  8.19 &  8.03 &  7.66 \\ 
		$10^{-10}$ &  2.34 &  4.86 &  5.95 &  7.60 &  7.85 &  8.27 &  8.27 &  8.32 &  8.24 &  8.18 \\ 
		$10^{-11}$ &  2.34 &  3.63 &  6.02 &  6.94 &  7.72 &  8.21 &  8.27 &  8.34 &  8.38 &  8.32 \\ 
		$10^{-12}$ &  2.34 &  3.63 &  6.02 &  6.94 &  7.72 &  8.22 &  8.28 &  8.36 &  8.41 &  8.40 \\ 
		$10^{-13}$ &  2.34 &  3.63 &  6.02 &  6.94 &  7.72 &  8.22 &  8.28 &  8.37 &  8.42 &  8.42 \\ 
		$10^{-14}$ &  2.34 &  5.10 &  6.78 &  7.83 &  8.27 &  8.40 &  8.44 &  8.45 &  8.45 &  8.45 \\ 
		\bottomrule
	\end{tabular}
	\caption{$\tnorm{u-u_N}_\eps \cdot 100 \big(N/(\KAK+1)\big)^2$ for $P_2$-FEM applied to
	Example~\ref{test:SunStynes} with certain $\eps$ and $\lambda = 0.25$.}
	\label{tab:KAKfem2}
	\end{center}
	\end{table}
	
	\begin{remark}
		In the numerically studied ranges ($\eps \in [10^{-14},1]$ and $\lambda \in \{0.005,0.25\}$)
		the $\eps$-dependent term in~\eqref{sigma} is dominant and thus $\KAK$ does not depend on $N$,
		cf. also Figure~\ref{fig:KAK}. Therefore, the logarithmic factor in the estimate of
		Theorem~\ref{th:energyErrorSunStynes} could not be seen in the numerical experiments.
		
		Moreover, the bounds on $\KAK$ of Section~\ref{sec:energySunStynes} suggest that in the
		tested parameter ranges $\KAK + 1$ logarithmically depends on $\eps$. But, as we have seen
		from the numerical studies, the energy norm seems to be too weak for the layers considered.
		This may explain why no logarithmic factor in $\eps$ is visible in the computational results.
		
		\begin{figure}
			\begin{center}
%
%
\begin{tikzpicture}[scale=0.9]

\begin{axis}[%
width=2.3in,
height=1.9in,
at={(0.809in,0.513in)},
scale only axis,
every outer x axis line/.append style={black},
every x tick label/.append style={font=\color{black}},
xmin=2,
xmax=12,
tick align=outside,
xlabel={$i$},
xmajorgrids,
every outer y axis line/.append style={black},
every y tick label/.append style={font=\color{black}},
ymin=0,
ymax=40,
ylabel={$j$},
ymajorgrids,
every outer z axis line/.append style={black},
every z tick label/.append style={font=\color{black}},
zmin=0,
zmax=20,
zmajorgrids,
view={-37.5}{30},
axis background/.style={fill=white},
title={$\KAK$ for $\lambda = 0.25$},
axis x line*=bottom,
axis y line*=left,
axis z line*=left
]

\addplot3[%
scatter, only marks,
shader=flat corner,draw=black,z buffer=sort,colormap/jet,mesh/rows=10]
table[row sep=crcr, point meta=\thisrow{c}] {%
x	y	z	c\\
3	0	1	1\\
3	1	1	1\\
3	2	1	1\\
3	3	2	2\\
3	4	2	2\\
3	5	3	3\\
3	6	3	3\\
3	7	4	4\\
3	8	4	4\\
3	9	5	5\\
3	10	5	5\\
3	11	5	5\\
3	12	5	5\\
3	13	5	5\\
3	14	5	5\\
3	15	5	5\\
3	16	5	5\\
3	17	5	5\\
3	18	5	5\\
3	19	5	5\\
3	20	5	5\\
3	21	5	5\\
3	22	5	5\\
3	23	5	5\\
3	24	5	5\\
3	25	5	5\\
3	26	5	5\\
3	27	5	5\\
3	28	5	5\\
3	29	5	5\\
3	30	5	5\\
3	31	5	5\\
3	32	5	5\\
3	33	5	5\\
3	34	5	5\\
3	35	5	5\\
3	36	5	5\\
3	37	5	5\\
3	38	5	5\\
3	39	5	5\\
3	40	5	5\\
3	41	5	5\\
3	42	5	5\\
3	43	5	5\\
3	44	5	5\\
3	45	5	5\\
3	46	5	5\\
3	47	5	5\\
3	48	5	5\\
3	49	5	5\\
3	50	5	5\\
4	0	1	1\\
4	1	1	1\\
4	2	1	1\\
4	3	2	2\\
4	4	2	2\\
4	5	3	3\\
4	6	3	3\\
4	7	4	4\\
4	8	4	4\\
4	9	5	5\\
4	10	5	5\\
4	11	6	6\\
4	12	6	6\\
4	13	6	6\\
4	14	7	7\\
4	15	7	7\\
4	16	7	7\\
4	17	7	7\\
4	18	7	7\\
4	19	7	7\\
4	20	7	7\\
4	21	7	7\\
4	22	7	7\\
4	23	7	7\\
4	24	7	7\\
4	25	7	7\\
4	26	7	7\\
4	27	7	7\\
4	28	7	7\\
4	29	7	7\\
4	30	7	7\\
4	31	7	7\\
4	32	7	7\\
4	33	7	7\\
4	34	7	7\\
4	35	7	7\\
4	36	7	7\\
4	37	7	7\\
4	38	7	7\\
4	39	7	7\\
4	40	7	7\\
4	41	7	7\\
4	42	7	7\\
4	43	7	7\\
4	44	7	7\\
4	45	7	7\\
4	46	7	7\\
4	47	7	7\\
4	48	7	7\\
4	49	7	7\\
4	50	7	7\\
5	0	1	1\\
5	1	1	1\\
5	2	1	1\\
5	3	2	2\\
5	4	2	2\\
5	5	3	3\\
5	6	3	3\\
5	7	4	4\\
5	8	4	4\\
5	9	5	5\\
5	10	5	5\\
5	11	6	6\\
5	12	6	6\\
5	13	6	6\\
5	14	7	7\\
5	15	7	7\\
5	16	8	8\\
5	17	8	8\\
5	18	8	8\\
5	19	8	8\\
5	20	8	8\\
5	21	8	8\\
5	22	8	8\\
5	23	8	8\\
5	24	8	8\\
5	25	8	8\\
5	26	8	8\\
5	27	8	8\\
5	28	8	8\\
5	29	8	8\\
5	30	8	8\\
5	31	8	8\\
5	32	8	8\\
5	33	8	8\\
5	34	8	8\\
5	35	8	8\\
5	36	8	8\\
5	37	8	8\\
5	38	8	8\\
5	39	8	8\\
5	40	8	8\\
5	41	8	8\\
5	42	8	8\\
5	43	8	8\\
5	44	8	8\\
5	45	8	8\\
5	46	8	8\\
5	47	8	8\\
5	48	8	8\\
5	49	8	8\\
5	50	8	8\\
6	0	1	1\\
6	1	1	1\\
6	2	1	1\\
6	3	2	2\\
6	4	2	2\\
6	5	3	3\\
6	6	3	3\\
6	7	4	4\\
6	8	4	4\\
6	9	5	5\\
6	10	5	5\\
6	11	6	6\\
6	12	6	6\\
6	13	6	6\\
6	14	7	7\\
6	15	7	7\\
6	16	8	8\\
6	17	8	8\\
6	18	9	9\\
6	19	9	9\\
6	20	10	10\\
6	21	10	10\\
6	22	10	10\\
6	23	10	10\\
6	24	10	10\\
6	25	10	10\\
6	26	10	10\\
6	27	10	10\\
6	28	10	10\\
6	29	10	10\\
6	30	10	10\\
6	31	10	10\\
6	32	10	10\\
6	33	10	10\\
6	34	10	10\\
6	35	10	10\\
6	36	10	10\\
6	37	10	10\\
6	38	10	10\\
6	39	10	10\\
6	40	10	10\\
6	41	10	10\\
6	42	10	10\\
6	43	10	10\\
6	44	10	10\\
6	45	10	10\\
6	46	10	10\\
6	47	10	10\\
6	48	10	10\\
6	49	10	10\\
6	50	10	10\\
7	0	1	1\\
7	1	1	1\\
7	2	1	1\\
7	3	2	2\\
7	4	2	2\\
7	5	3	3\\
7	6	3	3\\
7	7	4	4\\
7	8	4	4\\
7	9	5	5\\
7	10	5	5\\
7	11	6	6\\
7	12	6	6\\
7	13	6	6\\
7	14	7	7\\
7	15	7	7\\
7	16	8	8\\
7	17	8	8\\
7	18	9	9\\
7	19	9	9\\
7	20	10	10\\
7	21	10	10\\
7	22	11	11\\
7	23	11	11\\
7	24	11	11\\
7	25	11	11\\
7	26	11	11\\
7	27	11	11\\
7	28	11	11\\
7	29	11	11\\
7	30	11	11\\
7	31	11	11\\
7	32	11	11\\
7	33	11	11\\
7	34	11	11\\
7	35	11	11\\
7	36	11	11\\
7	37	11	11\\
7	38	11	11\\
7	39	11	11\\
7	40	11	11\\
7	41	11	11\\
7	42	11	11\\
7	43	11	11\\
7	44	11	11\\
7	45	11	11\\
7	46	11	11\\
7	47	11	11\\
7	48	11	11\\
7	49	11	11\\
7	50	11	11\\
8	0	1	1\\
8	1	1	1\\
8	2	1	1\\
8	3	2	2\\
8	4	2	2\\
8	5	3	3\\
8	6	3	3\\
8	7	4	4\\
8	8	4	4\\
8	9	5	5\\
8	10	5	5\\
8	11	6	6\\
8	12	6	6\\
8	13	6	6\\
8	14	7	7\\
8	15	7	7\\
8	16	8	8\\
8	17	8	8\\
8	18	9	9\\
8	19	9	9\\
8	20	10	10\\
8	21	10	10\\
8	22	11	11\\
8	23	11	11\\
8	24	11	11\\
8	25	12	12\\
8	26	12	12\\
8	27	13	13\\
8	28	13	13\\
8	29	13	13\\
8	30	13	13\\
8	31	13	13\\
8	32	13	13\\
8	33	13	13\\
8	34	13	13\\
8	35	13	13\\
8	36	13	13\\
8	37	13	13\\
8	38	13	13\\
8	39	13	13\\
8	40	13	13\\
8	41	13	13\\
8	42	13	13\\
8	43	13	13\\
8	44	13	13\\
8	45	13	13\\
8	46	13	13\\
8	47	13	13\\
8	48	13	13\\
8	49	13	13\\
8	50	13	13\\
9	0	1	1\\
9	1	1	1\\
9	2	1	1\\
9	3	2	2\\
9	4	2	2\\
9	5	3	3\\
9	6	3	3\\
9	7	4	4\\
9	8	4	4\\
9	9	5	5\\
9	10	5	5\\
9	11	6	6\\
9	12	6	6\\
9	13	6	6\\
9	14	7	7\\
9	15	7	7\\
9	16	8	8\\
9	17	8	8\\
9	18	9	9\\
9	19	9	9\\
9	20	10	10\\
9	21	10	10\\
9	22	11	11\\
9	23	11	11\\
9	24	11	11\\
9	25	12	12\\
9	26	12	12\\
9	27	13	13\\
9	28	13	13\\
9	29	14	14\\
9	30	14	14\\
9	31	14	14\\
9	32	14	14\\
9	33	14	14\\
9	34	14	14\\
9	35	14	14\\
9	36	14	14\\
9	37	14	14\\
9	38	14	14\\
9	39	14	14\\
9	40	14	14\\
9	41	14	14\\
9	42	14	14\\
9	43	14	14\\
9	44	14	14\\
9	45	14	14\\
9	46	14	14\\
9	47	14	14\\
9	48	14	14\\
9	49	14	14\\
9	50	14	14\\
10	0	1	1\\
10	1	1	1\\
10	2	1	1\\
10	3	2	2\\
10	4	2	2\\
10	5	3	3\\
10	6	3	3\\
10	7	4	4\\
10	8	4	4\\
10	9	5	5\\
10	10	5	5\\
10	11	6	6\\
10	12	6	6\\
10	13	6	6\\
10	14	7	7\\
10	15	7	7\\
10	16	8	8\\
10	17	8	8\\
10	18	9	9\\
10	19	9	9\\
10	20	10	10\\
10	21	10	10\\
10	22	11	11\\
10	23	11	11\\
10	24	11	11\\
10	25	12	12\\
10	26	12	12\\
10	27	13	13\\
10	28	13	13\\
10	29	14	14\\
10	30	14	14\\
10	31	15	15\\
10	32	15	15\\
10	33	16	16\\
10	34	16	16\\
10	35	16	16\\
10	36	16	16\\
10	37	16	16\\
10	38	16	16\\
10	39	16	16\\
10	40	16	16\\
10	41	16	16\\
10	42	16	16\\
10	43	16	16\\
10	44	16	16\\
10	45	16	16\\
10	46	16	16\\
10	47	16	16\\
10	48	16	16\\
10	49	16	16\\
10	50	16	16\\
11	0	1	1\\
11	1	1	1\\
11	2	1	1\\
11	3	2	2\\
11	4	2	2\\
11	5	3	3\\
11	6	3	3\\
11	7	4	4\\
11	8	4	4\\
11	9	5	5\\
11	10	5	5\\
11	11	6	6\\
11	12	6	6\\
11	13	6	6\\
11	14	7	7\\
11	15	7	7\\
11	16	8	8\\
11	17	8	8\\
11	18	9	9\\
11	19	9	9\\
11	20	10	10\\
11	21	10	10\\
11	22	11	11\\
11	23	11	11\\
11	24	11	11\\
11	25	12	12\\
11	26	12	12\\
11	27	13	13\\
11	28	13	13\\
11	29	14	14\\
11	30	14	14\\
11	31	15	15\\
11	32	15	15\\
11	33	16	16\\
11	34	16	16\\
11	35	17	17\\
11	36	17	17\\
11	37	17	17\\
11	38	17	17\\
11	39	17	17\\
11	40	17	17\\
11	41	17	17\\
11	42	17	17\\
11	43	17	17\\
11	44	17	17\\
11	45	17	17\\
11	46	17	17\\
11	47	17	17\\
11	48	17	17\\
11	49	17	17\\
11	50	17	17\\
12	0	1	1\\
12	1	1	1\\
12	2	1	1\\
12	3	2	2\\
12	4	2	2\\
12	5	3	3\\
12	6	3	3\\
12	7	4	4\\
12	8	4	4\\
12	9	5	5\\
12	10	5	5\\
12	11	6	6\\
12	12	6	6\\
12	13	6	6\\
12	14	7	7\\
12	15	7	7\\
12	16	8	8\\
12	17	8	8\\
12	18	9	9\\
12	19	9	9\\
12	20	10	10\\
12	21	10	10\\
12	22	11	11\\
12	23	11	11\\
12	24	11	11\\
12	25	12	12\\
12	26	12	12\\
12	27	13	13\\
12	28	13	13\\
12	29	14	14\\
12	30	14	14\\
12	31	15	15\\
12	32	15	15\\
12	33	16	16\\
12	34	16	16\\
12	35	17	17\\
12	36	17	17\\
12	37	17	17\\
12	38	18	18\\
12	39	18	18\\
12	40	19	19\\
12	41	19	19\\
12	42	19	19\\
12	43	19	19\\
12	44	19	19\\
12	45	19	19\\
12	46	19	19\\
12	47	19	19\\
12	48	19	19\\
12	49	19	19\\
12	50	19	19\\
};
\end{axis}
\end{tikzpicture}%
%
%
\begin{tikzpicture}[scale = 0.9]

\begin{axis}[%
width=2.3in,
height=1.9in,
at={(0.809in,0.513in)},
scale only axis,
every outer x axis line/.append style={black},
every x tick label/.append style={font=\color{black}},
xmin=2,
xmax=12,
tick align=outside,
xlabel={$i$},
xmajorgrids,
every outer y axis line/.append style={black},
every y tick label/.append style={font=\color{black}},
ymin=0,
ymax=40,
ylabel={$j$},
ymajorgrids,
every outer z axis line/.append style={black},
every z tick label/.append style={font=\color{black}},
zmin=0,
zmax=20,
zmajorgrids,
view={-37.5}{30},
axis background/.style={fill=white},
title={$\KAK$ for $\lambda = 0.005$},
axis x line*=bottom,
axis y line*=left,
axis z line*=left
]

\addplot3[%
scatter,only marks,
shader=flat corner,draw=black,z buffer=sort,colormap/jet,mesh/rows=10]
table[row sep=crcr, point meta=\thisrow{c}] {%
x	y	z	c\\
3	0	1	1\\
3	1	1	1\\
3	2	1	1\\
3	3	2	2\\
3	4	2	2\\
3	5	3	3\\
3	6	3	3\\
3	7	4	4\\
3	8	4	4\\
3	9	5	5\\
3	10	5	5\\
3	11	5	5\\
3	12	5	5\\
3	13	5	5\\
3	14	5	5\\
3	15	5	5\\
3	16	5	5\\
3	17	5	5\\
3	18	5	5\\
3	19	5	5\\
3	20	5	5\\
3	21	5	5\\
3	22	5	5\\
3	23	5	5\\
3	24	5	5\\
3	25	5	5\\
3	26	5	5\\
3	27	5	5\\
3	28	5	5\\
3	29	5	5\\
3	30	5	5\\
3	31	5	5\\
3	32	5	5\\
3	33	5	5\\
3	34	5	5\\
3	35	5	5\\
3	36	5	5\\
3	37	5	5\\
3	38	5	5\\
3	39	5	5\\
3	40	5	5\\
3	41	5	5\\
3	42	5	5\\
3	43	5	5\\
3	44	5	5\\
3	45	5	5\\
3	46	5	5\\
3	47	5	5\\
3	48	5	5\\
3	49	5	5\\
3	50	5	5\\
4	0	1	1\\
4	1	1	1\\
4	2	1	1\\
4	3	2	2\\
4	4	2	2\\
4	5	3	3\\
4	6	3	3\\
4	7	4	4\\
4	8	4	4\\
4	9	5	5\\
4	10	5	5\\
4	11	6	6\\
4	12	6	6\\
4	13	7	7\\
4	14	7	7\\
4	15	7	7\\
4	16	7	7\\
4	17	7	7\\
4	18	7	7\\
4	19	7	7\\
4	20	7	7\\
4	21	7	7\\
4	22	7	7\\
4	23	7	7\\
4	24	7	7\\
4	25	7	7\\
4	26	7	7\\
4	27	7	7\\
4	28	7	7\\
4	29	7	7\\
4	30	7	7\\
4	31	7	7\\
4	32	7	7\\
4	33	7	7\\
4	34	7	7\\
4	35	7	7\\
4	36	7	7\\
4	37	7	7\\
4	38	7	7\\
4	39	7	7\\
4	40	7	7\\
4	41	7	7\\
4	42	7	7\\
4	43	7	7\\
4	44	7	7\\
4	45	7	7\\
4	46	7	7\\
4	47	7	7\\
4	48	7	7\\
4	49	7	7\\
4	50	7	7\\
5	0	1	1\\
5	1	1	1\\
5	2	1	1\\
5	3	2	2\\
5	4	2	2\\
5	5	3	3\\
5	6	3	3\\
5	7	4	4\\
5	8	4	4\\
5	9	5	5\\
5	10	5	5\\
5	11	6	6\\
5	12	6	6\\
5	13	7	7\\
5	14	7	7\\
5	15	8	8\\
5	16	8	8\\
5	17	8	8\\
5	18	8	8\\
5	19	8	8\\
5	20	8	8\\
5	21	8	8\\
5	22	8	8\\
5	23	8	8\\
5	24	8	8\\
5	25	8	8\\
5	26	8	8\\
5	27	8	8\\
5	28	8	8\\
5	29	8	8\\
5	30	8	8\\
5	31	8	8\\
5	32	8	8\\
5	33	8	8\\
5	34	8	8\\
5	35	8	8\\
5	36	8	8\\
5	37	8	8\\
5	38	8	8\\
5	39	8	8\\
5	40	8	8\\
5	41	8	8\\
5	42	8	8\\
5	43	8	8\\
5	44	8	8\\
5	45	8	8\\
5	46	8	8\\
5	47	8	8\\
5	48	8	8\\
5	49	8	8\\
5	50	8	8\\
6	0	1	1\\
6	1	1	1\\
6	2	1	1\\
6	3	2	2\\
6	4	2	2\\
6	5	3	3\\
6	6	3	3\\
6	7	4	4\\
6	8	4	4\\
6	9	5	5\\
6	10	5	5\\
6	11	6	6\\
6	12	6	6\\
6	13	7	7\\
6	14	7	7\\
6	15	8	8\\
6	16	8	8\\
6	17	9	9\\
6	18	9	9\\
6	19	10	10\\
6	20	10	10\\
6	21	10	10\\
6	22	10	10\\
6	23	10	10\\
6	24	10	10\\
6	25	10	10\\
6	26	10	10\\
6	27	10	10\\
6	28	10	10\\
6	29	10	10\\
6	30	10	10\\
6	31	10	10\\
6	32	10	10\\
6	33	10	10\\
6	34	10	10\\
6	35	10	10\\
6	36	10	10\\
6	37	10	10\\
6	38	10	10\\
6	39	10	10\\
6	40	10	10\\
6	41	10	10\\
6	42	10	10\\
6	43	10	10\\
6	44	10	10\\
6	45	10	10\\
6	46	10	10\\
6	47	10	10\\
6	48	10	10\\
6	49	10	10\\
6	50	10	10\\
7	0	1	1\\
7	1	1	1\\
7	2	1	1\\
7	3	2	2\\
7	4	2	2\\
7	5	3	3\\
7	6	3	3\\
7	7	4	4\\
7	8	4	4\\
7	9	5	5\\
7	10	5	5\\
7	11	6	6\\
7	12	6	6\\
7	13	7	7\\
7	14	7	7\\
7	15	8	8\\
7	16	8	8\\
7	17	9	9\\
7	18	9	9\\
7	19	10	10\\
7	20	10	10\\
7	21	11	11\\
7	22	11	11\\
7	23	11	11\\
7	24	11	11\\
7	25	11	11\\
7	26	11	11\\
7	27	11	11\\
7	28	11	11\\
7	29	11	11\\
7	30	11	11\\
7	31	11	11\\
7	32	11	11\\
7	33	11	11\\
7	34	11	11\\
7	35	11	11\\
7	36	11	11\\
7	37	11	11\\
7	38	11	11\\
7	39	11	11\\
7	40	11	11\\
7	41	11	11\\
7	42	11	11\\
7	43	11	11\\
7	44	11	11\\
7	45	11	11\\
7	46	11	11\\
7	47	11	11\\
7	48	11	11\\
7	49	11	11\\
7	50	11	11\\
8	0	1	1\\
8	1	1	1\\
8	2	1	1\\
8	3	2	2\\
8	4	2	2\\
8	5	3	3\\
8	6	3	3\\
8	7	4	4\\
8	8	4	4\\
8	9	5	5\\
8	10	5	5\\
8	11	6	6\\
8	12	6	6\\
8	13	7	7\\
8	14	7	7\\
8	15	8	8\\
8	16	8	8\\
8	17	9	9\\
8	18	9	9\\
8	19	10	10\\
8	20	10	10\\
8	21	11	11\\
8	22	11	11\\
8	23	12	12\\
8	24	12	12\\
8	25	13	13\\
8	26	13	13\\
8	27	13	13\\
8	28	13	13\\
8	29	13	13\\
8	30	13	13\\
8	31	13	13\\
8	32	13	13\\
8	33	13	13\\
8	34	13	13\\
8	35	13	13\\
8	36	13	13\\
8	37	13	13\\
8	38	13	13\\
8	39	13	13\\
8	40	13	13\\
8	41	13	13\\
8	42	13	13\\
8	43	13	13\\
8	44	13	13\\
8	45	13	13\\
8	46	13	13\\
8	47	13	13\\
8	48	13	13\\
8	49	13	13\\
8	50	13	13\\
9	0	1	1\\
9	1	1	1\\
9	2	1	1\\
9	3	2	2\\
9	4	2	2\\
9	5	3	3\\
9	6	3	3\\
9	7	4	4\\
9	8	4	4\\
9	9	5	5\\
9	10	5	5\\
9	11	6	6\\
9	12	6	6\\
9	13	7	7\\
9	14	7	7\\
9	15	8	8\\
9	16	8	8\\
9	17	9	9\\
9	18	9	9\\
9	19	10	10\\
9	20	10	10\\
9	21	11	11\\
9	22	11	11\\
9	23	12	12\\
9	24	12	12\\
9	25	13	13\\
9	26	13	13\\
9	27	14	14\\
9	28	14	14\\
9	29	14	14\\
9	30	14	14\\
9	31	14	14\\
9	32	14	14\\
9	33	14	14\\
9	34	14	14\\
9	35	14	14\\
9	36	14	14\\
9	37	14	14\\
9	38	14	14\\
9	39	14	14\\
9	40	14	14\\
9	41	14	14\\
9	42	14	14\\
9	43	14	14\\
9	44	14	14\\
9	45	14	14\\
9	46	14	14\\
9	47	14	14\\
9	48	14	14\\
9	49	14	14\\
9	50	14	14\\
10	0	1	1\\
10	1	1	1\\
10	2	1	1\\
10	3	2	2\\
10	4	2	2\\
10	5	3	3\\
10	6	3	3\\
10	7	4	4\\
10	8	4	4\\
10	9	5	5\\
10	10	5	5\\
10	11	6	6\\
10	12	6	6\\
10	13	7	7\\
10	14	7	7\\
10	15	8	8\\
10	16	8	8\\
10	17	9	9\\
10	18	9	9\\
10	19	10	10\\
10	20	10	10\\
10	21	11	11\\
10	22	11	11\\
10	23	12	12\\
10	24	12	12\\
10	25	13	13\\
10	26	13	13\\
10	27	14	14\\
10	28	14	14\\
10	29	15	15\\
10	30	15	15\\
10	31	16	16\\
10	32	16	16\\
10	33	16	16\\
10	34	16	16\\
10	35	16	16\\
10	36	16	16\\
10	37	16	16\\
10	38	16	16\\
10	39	16	16\\
10	40	16	16\\
10	41	16	16\\
10	42	16	16\\
10	43	16	16\\
10	44	16	16\\
10	45	16	16\\
10	46	16	16\\
10	47	16	16\\
10	48	16	16\\
10	49	16	16\\
10	50	16	16\\
11	0	1	1\\
11	1	1	1\\
11	2	1	1\\
11	3	2	2\\
11	4	2	2\\
11	5	3	3\\
11	6	3	3\\
11	7	4	4\\
11	8	4	4\\
11	9	5	5\\
11	10	5	5\\
11	11	6	6\\
11	12	6	6\\
11	13	7	7\\
11	14	7	7\\
11	15	8	8\\
11	16	8	8\\
11	17	9	9\\
11	18	9	9\\
11	19	10	10\\
11	20	10	10\\
11	21	11	11\\
11	22	11	11\\
11	23	12	12\\
11	24	12	12\\
11	25	13	13\\
11	26	13	13\\
11	27	14	14\\
11	28	14	14\\
11	29	15	15\\
11	30	15	15\\
11	31	16	16\\
11	32	16	16\\
11	33	17	17\\
11	34	17	17\\
11	35	17	17\\
11	36	17	17\\
11	37	17	17\\
11	38	17	17\\
11	39	17	17\\
11	40	17	17\\
11	41	17	17\\
11	42	17	17\\
11	43	17	17\\
11	44	17	17\\
11	45	17	17\\
11	46	17	17\\
11	47	17	17\\
11	48	17	17\\
11	49	17	17\\
11	50	17	17\\
12	0	1	1\\
12	1	1	1\\
12	2	1	1\\
12	3	2	2\\
12	4	2	2\\
12	5	3	3\\
12	6	3	3\\
12	7	4	4\\
12	8	4	4\\
12	9	5	5\\
12	10	5	5\\
12	11	6	6\\
12	12	6	6\\
12	13	7	7\\
12	14	7	7\\
12	15	8	8\\
12	16	8	8\\
12	17	9	9\\
12	18	9	9\\
12	19	10	10\\
12	20	10	10\\
12	21	11	11\\
12	22	11	11\\
12	23	12	12\\
12	24	12	12\\
12	25	13	13\\
12	26	13	13\\
12	27	14	14\\
12	28	14	14\\
12	29	15	15\\
12	30	15	15\\
12	31	16	16\\
12	32	16	16\\
12	33	17	17\\
12	34	17	17\\
12	35	18	18\\
12	36	18	18\\
12	37	19	19\\
12	38	19	19\\
12	39	19	19\\
12	40	19	19\\
12	41	19	19\\
12	42	19	19\\
12	43	19	19\\
12	44	19	19\\
12	45	19	19\\
12	46	19	19\\
12	47	19	19\\
12	48	19	19\\
12	49	19	19\\
12	50	19	19\\
};
\end{axis}
\end{tikzpicture}%
				\caption{$\KAK$ for $N = 2^i$, $\eps = 10^{-j}$, $k=2$, and $\lambda= 0.25$ (left)
				and $\lambda = 0.005$ (right).}
				\label{fig:KAK}
			\end{center}
		\end{figure}
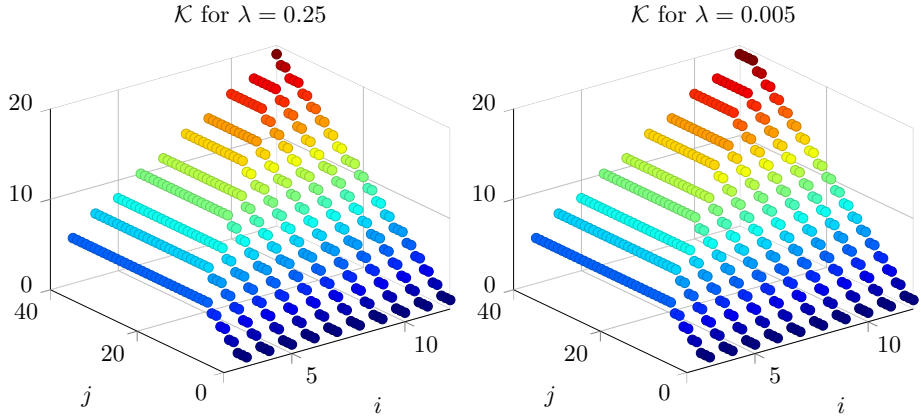
	\end{remark}
		
	In summary, our numerical experiments confirm the theoretical results of
	Section~\ref{sec:energySunStynes} and Section~\ref{sec:SDFEM}. However,
	the calculations suggest that the given bound for the $L^2$-norm error
	is not optimal yet. Also some interesting effects differing between 
	odd and even order elements could be seen when $\lambda$ is not too
	small. Furthermore, the computational results indicate that the numerical
	and theoretical study of other (balanced) norms for the layers
	of ``cusp''-type would be very interesting and should be object of further
	research.

	\section*{Acknowledgement}
	\noindent The author would like to thank Hans-G{\"o}rg Roos and Sebastian Franz for helpful
	comments and discussions.
	
	\bibliographystyle{plain}
	\bibliography{sunStynesBiblio}
	
\end{document}